\documentclass{m2an}
\usepackage{amssymb}
\usepackage{amsmath}
\usepackage{graphicx}
\usepackage{latexsym}
\usepackage{wrapfig}
\usepackage{placeins}
\usepackage{float}
\usepackage{comment}
\usepackage{caption}
\usepackage{subcaption}
\usepackage{cleveref}
\usepackage[dvipsnames]{xcolor}

\renewcommand{\div}{\operatorname{div}}

\newcommand{\RT}{\mathcal{RT}}

\newcommand{\jump}[1]{[#1]}

\newtheorem{remark}{Remark}

\renewcommand{\div}{{\hbox{div}}}

\renewcommand{\c}{\frac{1}{\delta x^2}}

\newcommand{\dx}{\, dx}
\newcommand{\ds}{\, ds}

\newcommand{\T}{\mathcal{T}_h}
\newcommand{\C}{\mathcal{C}_h}
\newcommand{\D}{\mathcal{D}_h}
\newcommand{\M}{\mathcal{M}_h}
\newcommand{\F}{\mathcal{F}_h}
\newcommand{\N}{\mathcal{N}_h}

\newtheorem{assumption}{Assumption}
\renewcommand{\div}{{\hbox{div}}}

\renewcommand{\c}{\frac{1}{\delta x^2}}
\newcommand{\sigIRT}{{\sigma}_h}
\newcommand{\IRT}{\mathcal{IRT}}

\newcommand{\Fg}{\mathcal{F}_g}
\newcommand{\jF}[1]{[\![#1]\!]}
 
\newcommand{\nrmL}[2][1]{\|#2\|_{#1}}

\newcommand{\Tp}{{T'}}

\newtheorem{example}{Example}[section]
\newcommand{\Blue}[1]{\textcolor{black}{#1}}
\newcommand{\Red}[1]{\textcolor{black}{#1}}
\newcommand{\Magenta}[1]{\textcolor{black}{#1}}
\begin{document}
%%-----------------------------
%%      the top matter
%%-----------------------------
\title{Elliptic interface problem approximated by CutFEM: II. A posteriori error analysis based on equilibrated fluxes}

\author{Daniela Capatina}\address{LMAP CNRS UMR 5142, University of Pau, 64013 Pau, France}
\author{Aimene Gouasmi}\address{LMAP CNRS UMR 5142, University of Pau, 64013 Pau, France}
\date{December 4th, 2025}
\begin{abstract} 
This paper investigates an elliptic interface problem with discontinuous diffusion coefficients \Red{discretized by finite elements} on unfitted meshes, employing the CutFEM method. The main contribution is the a posteriori error analysis based on equilibrated fluxes belonging to the immersed Raviart-Thomas space. We establish sharp reliability and local efficiency of a new error estimator, which includes both volume and interface terms, carefully tracking the dependence of the efficiency constant on the diffusion coefficients and the mesh/interface configuration. Numerical results highlight the robustness of the approach.
\end{abstract}

\begin{resume} On consid\`ere un probl\`eme elliptique d'interface avec des coefficients discontinus, \Red{discr\'etis\'e par \'el\'ements finis} sur des maillages non-conformes en utilisant la m\'ethode CutFEM. La principale contribution est l'analyse d'erreur a posteriori bas\'ee sur des flux \'equlibr\'es appartenant \`a l'espace de Raviart-Thomas immerg\'e. Nous \'etablissons la fiabilit\'e et l'efficacit\'e locale d'un nouvel estimateur d'erreur, qui inclut \`a la fois des termes sur les cellules du maillage et sur l'interface, en explicitant la d\'ependance des constantes d'efficacit\'e par rapport aux coefficients de diffusion et \`a la configuration maillage/interface. Des r\'esultats num\'eriques illustrent la robustesse de l'approche propos\'ee. \end{resume}

\maketitle
%%-----------------------------
%%      your text
%%-----------------------------
\section*{Introduction}

The study of a posteriori error estimation \Red{for finite element methods} has been an active area of research for several decades, as demonstrated by the extensive literature on the topic (see, for example, \cite{verfurth1994posteriori, ainsworth1997posteriori,  ern2010guaranteed} and many others). This technique is particularly valuable in adaptive mesh refinement procedures, commonly employed for problems with singularities, discontinuities, or sharp derivatives. The a posteriori analysis quantifies the error incurred and provides guidance on how to improve the accuracy of the solutions \cite{verfurth1994posteriori}. Several types of a posteriori error estimators exist, such as those based on the residuals of the governing differential equations or those based on equilibrated fluxes, which generally yield sharp reliability bounds. Various flux-based a posteriori estimators for different finite element methods have been studied by many researchers, see e.g.  \cite{LaLe:83, oden1989toward,  AiOd:93, ainsworth1997posteriori,  DeMe:99, braess2008equilibrated, Ve:09, vohralik2011guaranteed, CaZh:11, ErnVo2015, Dana2016}. \Blue{Using a partition of unity, the pioneering work of \cite{LaLe:83} introduced a local strategy that reduces the construction of equilibrated fluxes to vertex-patch computations. For the $P^1$-continuous finite element approximation of the Poisson equation, an equilibrated flux in the lowest-order Raviart–Thomas space was explicitly built in \cite{braess2008equilibrated, DeMe:99} and  via local corrections in the broken Raviart–Thomas space. However, this approach does not yield a robust equilibrated estimator with respect to coefficient jumps unless a constrained minimization is introduced; see \cite{CaZh:11}. In \cite{ErnVo2015}, a unified partition-of-unity method was developed, requiring the solution of local mixed problems on each vertex patch. As an exception to partition-of-unity-based approaches, we refer to \cite{Dana2016} for a Poisson problem.}

We consider here a finite element approximation of an elliptic interface problem on unfitted meshes, where the interface $\Gamma$ is taken into account by the CutFEM method \cite{CutFEM2, CutFEM}. In \cite{Article1}, we have defined a locally conservative flux belonging to an immersed Raviart-Thomas space (cf. \cite{IRT}), \Red{following the approach developed in \cite{Aimene} for fitted meshes. We have also introduced} a global error estimator consisting of a standard term $\eta$ -the weighted $L^2$-norm of the difference between the
equilibrated and numerical fluxes- plus a new interface term $\eta_\Gamma$, accounting for discontinuities of the discrete solution across the interface and of the equilibrated flux across the cut edges. In this paper, we provide a detailed a posteriori error analysis based on this estimator, establishing its robust reliability and local efficiency. To the best of our knowledge, these topics are largely unexplored in the context of CutFEM for interface problems. \Red{An a posteriori error analysis based on conservative fluxes was carried out in \cite{DanaCuiyu} for CutFEM solutions of a Poisson boundary problem on unfitted meshes, but the challenges and techniques required for interface problems are more involved.}

We first establish the robust reliability and bound the weighted $H^1$ semi-norm of the error by the global estimator $\eta +\eta_{\Gamma}$. The reliability constant of $\eta$ is equal to $1$, in agreement with well-known results for equilibrated flux-based estimators, while the constant in front of $\eta_{\Gamma}$ and the higher-order term is independent of the mesh size, the diffusion coefficients and the mesh/interface geometry. The proof requires introducing a continuous approximation of the CutFEM solution such that the corresponding interpolation error is bounded by the interface estimator $\eta_{\Gamma}$. 

Then we focus on the proof of the local efficiency for each of the three local error indicators, aiming to obtain the best possible constants with respect to the various parameters. The treatment of the cut cells is challenging, since the difference between the equilibrated and the numerical fluxes is only piecewise Raviart-Thomas. Therefore, we first investigate this space and construct specific basis functions, in relation \Red{to} the immersed Raviart-Thomas space \cite{IRT}. Their estimates allow us to bound the flux-based indicators by the residual ones; by adapting the Verf\"urth argument \cite{verfurth1994posteriori} to the cut cells, we are next able to establish the efficiency, with constants that depend explicitly on the ratio of the diffusion coefficients and the mesh/interface geometry. 

Finally, we carry out some numerical tests. The adaptive mesh refinement is validated in \cite{Article1}, so we focus here on illustrating the estimator's behavior. We numerically observe that $\eta_{\Gamma}$ has a minor influence and thus, can be neglected in the mesh refinement, and that the estimator is robust with respect to the various parameters.

The paper is organized as follows. In Section \ref{sec:pb_modele}, we give the continuous and discrete formulations, while in Section \ref{sec: A posteriori_estim}, the equilibrated flux and the a posteriori error estimator introduced in \cite{Article1} are recalled. Section \ref{sec:reliab} is devoted to the proof of the reliability bound. In the next two sections, we establish some results useful to prove the efficiency: Section  \ref{sec:decomp_RT} is dedicated to the properties of a piecewise Raviart-Thomas space on the cut cells, while in Section \ref{sec:bound_tau} the decomposition of the flux correction in the previous space is obtained. Then we establish the local efficiency in Section \ref{sec:eff}. Finally, several numerical experiments are presented in Section \ref{sec:num_sim}. 
%The paper ends with Appendix \ref{sec:appendixB}, which deals with a technical point in the efficiency proof.

\section{Continuous and discrete problems}\label{sec:pb_modele}
The model problem is given by: 
\begin{equation}\label{intro problem_2}
\left \{
\begin{array}{ll}
-\mathrm{div}(K\nabla u_i)=f & \text{ in } \Omega^i \quad (i=1, 2)\Red{,}\\
u= 0 &\text{ on } \partial \Omega\Red{,} \\
\left[u\right]=0,\,\, \left[K\nabla u\cdot n_\Gamma\right] = g & \text{ on }\Gamma\Red{.}
\end{array}
\right.
\end{equation}
where $\Omega$ is a 2D polygonal domain and $\Gamma$ a sufficiently smooth interface, separating $\Omega$ into two disjoint subdomains $\Omega^{1}$ and $\Omega^{2}$. We denote by $n_{\Gamma}$ the unit normal vector to $\Gamma$ oriented from $\Omega^1$ to $\Omega^2$; the jumps across $\Gamma$ are defined by $[u]=u_1-u_2$,
where $u_{|\Omega^i}=u_i$. We suppose $f\in L^2(\Omega)$, $g \in L^2(\Gamma)$ and, for the sake of simplicity, we assume $K_{|\Omega^i}=K_i=k_i I$ with $k_i>0$, for $i=1,2$. 
We consider the following well-posed variational formulation: Find $u\in H^1_0(\Omega) $ such that
\begin{equation}\label{eq: continuous_problem}
\int_{\Omega}K\nabla u \cdot \nabla v\dx= \int_{\Omega}f v\dx+\int_{\Gamma}gv\ds,\qquad \forall v\in H^1_0(\Omega).
\end{equation}  

In view of the numerical approximation of \eqref{eq: continuous_problem} by the CutFEM method \cite{CutFEM2}, we introduce some notation (see also \cite{Article1}). Let $\mathcal{T}_h$ \Red{be} a regular triangular mesh of $\Omega$, whose elements are closed sets, $\mathcal{F}_h$ the set of edges, $h_T$ and $h_F$ the diameter of $T \in \mathcal{T}_h$ and $F \in \mathcal{F}_h$, respectively. \Blue{The shape regularity of the triangles is implicitly used in several places in the manuscript.} For an interior edge $F$, $n_F$ is a fixed unit normal vector, oriented from $T_F^-$ to $T_F^+$ (the two triangles that share $F$). If $F \subset \partial \Omega$, then $n_F$ is the outward normal vector to $\Omega$ while if $F \subset \Gamma$, then $n_F = n_{\Gamma}$.  For $i = 1, 2$, we define:
\begin{equation*}
	\mathcal{T}_h^{i} = \big\{ T \in \mathcal{T}_h \ ; \ T \cap \Omega^{i} \neq \emptyset \big\}, \quad \mathcal{F}_h^{i} = \big\{ F \in \mathcal{F}_h \ ; \ F \cap \Omega^{i} \neq \emptyset \big\},\quad \Omega_h^i = \bigcup_{T \in \mathcal{T}_h^i} T.
\end{equation*}
For $T\in\T$ and $F\in\F$, let \Red{us define} the patches $\displaystyle\Delta_T=\bigcup_{N\in \mathcal N_T}\omega_N$ and $\displaystyle\Delta_F=\bigcup_{T,\,\partial T\supset F}\Delta_T$, where $\mathcal N_T$ denotes the set of vertices of $T$ and $\omega_N$ the set of triangles sharing the node $N$. As regards the cut elements, let:
\begin{equation*}
\begin{split}
\mathcal{T}_h^{\Gamma} = \big\{ T \in \mathcal{T}_h \ ; \ \overset{\Blue{\circ}}T \cap \Gamma \neq \emptyset \big\}, &\quad \mathcal{F}_h^{\Gamma} = \big\{ F \in \mathcal{F}_h \ ; \ \overset{\Blue{\circ}}F \cap \Gamma \neq \emptyset \big\}, \\
T^i = T \cap \Omega^i\,\, \forall T \in \mathcal{T}_h^{\Gamma},&\quad F^i = F \cap \Omega^i \,\, \forall F \in \mathcal{F}_h^{\Gamma}.
\end{split}
\end{equation*}
\Blue{The triangles $T$ of $\mathcal T_h$ such that $T\cap \Gamma$ is a single vertex or a whole side of $T$ are treated like non-cut elements.} 

\Blue{For technical reasons related to the unisolvency of the immersed Raviart-Thomas space, we assume, as in \cite{IRT}, that all cut triangles are acute. This assumption is not overly restrictive, since the cut elements are likely to be refined and the first refinement step transforms an obtuse triangle into an acute one (because we use the longest edge refinement method).}

\Blue{In order to focus on the a posteriori error analysis based on the new flux, we neglect in this paper the interface approximation error and} we assume that $\Gamma$ is a polygonal line, such that for each $T \in \mathcal{T}_h^{\Gamma}$, the intersection $\Gamma_T := T \cap \Gamma$ is a \Blue{segment. We do not allow for interface corners inside an element; in such a situation, the triangle should be split into two sub-triangles, so that $\Gamma_T$ becomes a segment in each one.}

For a function $v$, sufficiently smooth on each $\Omega^i$ but discontinuous across $\Gamma$, we denote by $v_i := v_{|\Omega^i}$ and define (cf. \cite{Ern}) the two following means on $ \Gamma$:
\begin{equation*}
\{v\} = \omega_1 v_1 + \omega_2 v_2, \quad \{v\}^* = \omega_2 v_1 + \omega_1 v_2\quad \text{with}\quad \omega_1 = \frac{k_2}{k_2 + k_1}, \quad \omega_2 = \frac{k_1}{k_1 + k_2}.
\end{equation*}
It is also useful to introduce the harmonic mean $k_{\Gamma}=k_1 k_2/(k_1 + k_2)$. Furthermore, we introduce the arithmetic mean and the jump across an interior edge $F$ as follows:
$$\langle v \rangle = \frac{1}{2} (v^- + v ^+), \quad [\! [v]\!] = v^- - v^+,\quad [\! [\partial_n v]\!] =[\! [\nabla v]\!] \cdot n_F.$$
For boundary edges, we set $\langle v \rangle = [\![v]\!] = v$. We use the symbols $\lesssim$ and $\approx$ to indicate the existence of generic constants that are independent of the mesh size, the interface geometry and the diffusion coefficients. Finally, we recall the following trace inequality on a cut element $T\in \T^{\Gamma}$:
\begin{equation}\label{eq: trace_in_cut}
\forall v\in H^1(T),\quad \nrmL[\Gamma_{T}]{v}^2\lesssim\bigg( h_{T}^{-1}\nrmL[T]{v}^2+ h_{ T}\nrmL[T]{\nabla v}^2\bigg).
\end{equation}

The CutFEM approximation of \eqref{eq: continuous_problem} reads: Find $u_h=(u_{h,1},u_{h,2}) \in \C^1\times \C^2 $ such that
\begin{equation}\label{eq: Cutfem_Formulation}
	a_h(u_h,v_h)=l_h(v_h) \quad \forall v_h\in \C^1\times \C^2,
\end{equation}
where $ \C^i=\big\{ v \in H^1(\Omega^i_h);\,v_{|(\partial \Omega^i\setminus\Gamma)}=0,\, \, v_{|T}\in P^1(T),\,\, \forall \ T\in  \mathcal{T}_h^{i} \big\}$ for $i=1,2$ and 
\begin{equation*}
\begin{split}
    a_h(u_h,v_h)=& \sum_{i=1}^{2}\bigg(\sum_{T\in\T^i}\int_{T^i}k_i\nabla u_{h,i} \cdot \nabla v_{h,i}\dx+\gamma_g j_i(u_{h,i},v_{h,i})\bigg)+a_{\Gamma}(u_h,v_h),\\
j_i(u_{h,i},v_{h,i})=&\sum_{F\in\mathcal{F}_g^i} h_F \int_{F}k_i[\![\partial_n u_{h,i}]\!][\![\partial_n v_{h,i}]\!]\ds\qquad (i=1,\,2),\\
a_{\Gamma}(u_h,v_h)=& \sum_{T\in\T^{\Gamma}}\int_{\Gamma_T}\bigg ( \frac{\gamma k_{\Gamma}}{h_T}[u_h][v_h]-\{K\nabla u_h\cdot n_{\Gamma}\}[v_h]-\{K\nabla v_h\cdot n_{\Gamma}\}[u_h]\bigg)\ds,\\
l_h(v_h)=& \sum_{i=1}^{2}\int_{\Omega^i}f v_{h,i}\dx+\sum_{T\in\mathcal{T}_h^{\Gamma}}\int_{\Gamma_T}g\{v_h\}^*\ds.
\end{split}
\end{equation*}
Here above, $\gamma >0$ and $\gamma_g>0$ are stabilization parameters, independent of the mesh, the interface and the diffusion coefficients, while $\mathcal{F}_g^i = \big\{ F \in \mathcal{F}_h^i \ \Blue{;} \ (T_F^+ \cup T_F^-) \cap \Gamma \neq \emptyset \big\}$.

It is well known that for $\gamma$ large enough, the discrete problem (\ref{eq: Cutfem_Formulation}) is well-posed.

It is useful to introduce the following norms, for any $v_h\in \C^1\times \C^2$:
\begin{equation*}
|v_h|_{1,K,h}^2= \sum_{i=1}^{2}\|k_i^{1/2}\nabla v_{h,i}\|^2_{\Omega^i},\,\,
\|v_h\|_h^2= |v_h|_{1,K,h}^2+\sum_{i=1}^{2}j_i(v_{h,i},v_{h,i})
+\sum_{T\in \mathcal{T}_h^{\Gamma}}\frac{k_{\Gamma}}{h_T}\|[v_h]\|^2_{\Gamma_T}
\end{equation*}
\Blue{where $\|\cdot\|_\omega$ is the $L^2$-norm on $\omega \subset \mathbb{R}^d$ ($1\le d\le 2$).} Let the local semi-norms on $T\in \T$ and $\Delta_T$:
$$|v_h|_{1,K,T}^2=  \sum_{i=1}^{2}\|k_i^{1/2}\nabla v_{h,i}\|^2_{T^i}, \qquad
j_{i,\Delta_T}(v_{h,i},v_{h,i}):=\Blue{\sum_{T'\in \Delta_T}}\sum_{F\in \mathcal F_g^i\cap \Blue{\mathcal F_{T'}}} h_F k_i\int_F \jF{\partial_n v_{h,i}}^2\ds,$$
\Blue{which can also be defined on any subset $\omega\subset \Omega$ that is a union of cells. Here above, $\mathcal F_T$ is the set of sides of $T$.}

In \cite{Article1}, an equivalent hybrid mixed formulation of \eqref{eq: Cutfem_Formulation} was proposed and studied. We recall below the definition of its Lagrange multiplier $\theta_h=(\theta_{h,1},\theta_{h,2})$, which was further used to reconstruct an equilibrated flux. For this purpose, let the spaces:
\begin{equation*}
\begin{split}
\D^i=&\{v\in L^2(\T^i);\, {v}_{|T}\in P^1(T)\ \forall T\in \T^i\}, \qquad (i=1,2),\\
\M^i=&\bigg\{\mu\in L^2(\F^i);\, {\mu}_{|F}\in P^1(F)\,\, \forall F\in \F ^i,\, 
\displaystyle\sum_{F\in\mathcal{F}_N}\mathfrak{s}_N^Fh_F{\mu}_{|F}(N)=0 \ \ \,\forall N\in \overset{\circ}{\mathcal {N}_h^i}\bigg \},
\end{split}     
\end{equation*}
where $\overset{\circ}{\mathcal {N}_h^i}$ denotes the set of nodes interior to $\Omega_h^i$, $\mathcal{F}_N$ the set of edges sharing the node $N$, and $\mathfrak{s}_N^F= \text{sign}(n_F,N)$ is equal to $1 (-1)$ if the orientation of $n_F$ with respect to the node $N$ is in clockwise (counter-clockwise) rotation. Let also the bilinear form
\begin{equation*}
b_h(\mu_h,v_h)=\displaystyle \sum_{i=1}^{2}\sum_{F\in\mathcal{F}_h^{i}} \frac{k_i h_F}{2}\sum_{N\in\mathcal{N}_F}{\mu_{h,i}}_{|F}(N)[\![v_{h,i}]\!](N),\,\, \forall \mu_h\in \M^1\times \M^2,\,v_h\in \D^1\times \D^2
\end{equation*}
where $\mathcal N_F$ is the set of vertices of the edge $F$. 

Then we consider the following problem: Find $\theta_h=(\theta_{h,1},\theta_{h,2}) \in \M^1\times \M^2$ s.t. 
\begin{equation*}
b_h(\theta_h,v_h)=l_h(v_h)-a_h(u_h,v_h)+\sum_{i=1}^{2} \sum_{F\in\mathcal{F}_h^{i}} \int_{F^i} \langle k_i\nabla u_{h,i}\cdot n_F\rangle [\![v_{h,i}]\!] \ds\quad \forall v_h \in \D^1\times \D^2.
\end{equation*}

It was shown in \cite{Article1} that it has a unique solution, which can be computed locally as sum of contributions $\theta_N^{i} $ defined on the patches $\omega_N^i=\omega_N\cap \T^i$ associated to the nodes $N$ of $\T^i$ : $ \theta_{h,i}=\displaystyle\sum_{N\in \N^i}\theta_N^{i}$ for $ 1\le i\le 2$. \Blue{More precisely, $\theta_{h,i}$ belongs to $\mathcal M_h^i$, lives on the sides of $\mathcal F_N\cap \mathcal F_h^i$ and is the unique solution of a small linear system associated to $N$, see \cite{Article1} for details.} In addition, we have the next bound \Blue{(see Theorem 2.6 of \cite{Article1})}:
\begin{lmm}\label{thm:bound_theta}
For $i=1,2$ and any node $N$ of $\T^i$, we have that:
\begin{equation*}
\begin{split}
\bigg(\sum_{F\in \mathcal F_N \cap \F^i} h_F k_i^{2}\|\theta_N^i\|_F^2\bigg)^{1/2} &\lesssim \sum_{T\in\omega_N^i} \sqrt{h_{T}}k_i\bigg(\| \jF{\partial_n u_{h,i}}\|_{\partial T\setminus \partial \omega_N^i} +\| \jF{\partial_n u_{h,i}}\|_{\partial T\cap \Fg^i}\bigg)\\
 &+\sum_{T\in\T^\Gamma\cap \omega_N^i}\left(\frac{{k_\Gamma}}{\sqrt{ h_T}}\|\jump{u_h}\|_{\Gamma_T}+ \sqrt{h_{T}}\|g-[K\nabla u_h\cdot n_\Gamma]\|_{\Gamma_T}\right)
 +\sum_{T\in\omega_N^i}h_T\|f\|_{T^i}.
\end{split} 
\end{equation*}
\end{lmm}

\section{Definition of the flux and the a posteriori error estimators}\label{sec: A posteriori_estim}
For the sake of simplicity, we assume from now on $g=0$ but the study can be extended to non-homogeneous transmission conditions. 
\Red{We recall that the lowest-order Raviart-Thomas space is defined as $\RT^0(\mathbb R^2)=\{(ax_1+b,ax_2+c);\,a,b,c\in \mathbb R\}$, where $x=(x_1,x_2)\in \mathbb R^2$.} \Blue{For $\omega\subset \mathbb R^2$, we set $\RT^0(\omega)=\{\psi_{|\omega};\, \psi\in \RT^0(\mathbb R^2)\}$.}

In \cite{Article1}, an equilibrated conservative flux $\sigma_h$ was defined in the immersed Raviart-Thomas space $\IRT^0(\mathcal T_h)$ of \cite{IRT}. We recall that on a cut cell $T$, the elements of this immersed space are piecewise $\RT^0$- functions $\psi$ satisfying:
\begin{equation*}
[\psi\cdot n_\Gamma] =0,\quad 
[K^{-1}\psi\cdot t_\Gamma](x_\Gamma)=0,\quad \div (\psi_{|T^1})= \div (\psi_{|T^2}),
\end{equation*}
where $x_\Gamma$ is an arbitrary point of $\Gamma_T$. On any $T\in \T$, the local degrees of freedom are the same as for the standard $\RT^0$-space, that is:   
\begin{equation}\label{eq:def_ddl_RT}
    N_{T,j}(\psi)=\frac{1}{|F_j|}\int_{F_j}\psi\cdot n_{T}\ds\quad   (1\le j\le 3),
\end{equation}
where $(F_j)_{1\le j\le 3}$ denote the edges of $T$ \Red{and $n_T$ is the outward unit normal vector to $T$}. The recovered flux $\sigma_h\in \IRT^0(\mathcal T_h)$ of \cite{Article1} is defined by imposing its degrees of freedom as follows:
\begin{equation*}
\begin{split}
\sigma_h\cdot n_F=&  \langle k_i\nabla u_{h,i}\cdot n_F\rangle - k_i\pi_F^0\theta_{h,i}, \quad \forall F\in \F^i\backslash\F^{\Gamma}\quad(i=1,2),\\
\int_F\sigma_h\cdot n_F\ds= &\sum_{i=1}^{2} \bigg (\int_{F^i} \langle k_i\nabla u_{h,i}\cdot n_F\rangle\ds-\Blue{\int_F} k_i \theta_{h,i}\Blue{\ds}\bigg), \quad \forall F\in \F^\Gamma,
\end{split}
\end{equation*}
\Blue{where $\pi^0_F$ denotes the $L^2(F)$-orthogonal projection on $P^0(F)$.} \Red{The definition of the degrees of freedom of the immersed Raviart-Thomas space $\mathcal {IRT}^O$ yields that $\displaystyle \int_F [\![\sigma_h\cdot n_F]\!]\ds=0$ on any interior side $F$, that is $\pi_0^F[\![\sigma_h\cdot n_F]\!]=0$.}
\Magenta{It was proved in Theorem 3.1 of \cite{Article1} that $\operatorname{div} \sigma_{h}=-f_{h}$, with $(f_h)_{|T}=\pi^0_T f$ for any $T\in\T$ and $ \pi^0_T $ the $L^2(T)$-orthogonal projection on $P^0(T)$.} 

We next set $\tau_h=K^{-1/2}(\sigma_h-K\nabla_h u_h)$ and define the local error estimator:
\begin{equation*}
    \eta_{T}
    =\|K^{-1/2}(\sigma_h - K\nabla_h u_h)\|_{T}=\|\tau_h \|_{T},\quad \forall T\in \mathcal T_h.
\end{equation*}
Since $u_h$ is discontinuous across $\Gamma$, we use the discrete gradient $\nabla_h$ in a cut triangle $T\in \T^\Gamma$; thus, $ \nabla_hu_h\in L^2(T)$ is defined by $
(\nabla_h {u_{h}})_{|T\cap\Omega^{i}} = (\nabla {u_{h,i}})_{|T\cap \Omega^{i}}$ for $i=1, 2$.

In addition, on the cut cells $T\in \T^\Gamma$ and cut edges $F\in \mathcal F_h^{\Gamma}$ we also consider
\begin{equation*}
    \tilde{\eta}_T=  \frac{\sqrt {h_T k_\Gamma}}{\sqrt {h_T^{min}|\Gamma_T|}}\| [u_h]\|_{\Gamma_T},\qquad \eta_F=\frac{\sqrt{h_F}}{\sqrt{k_{\Gamma}}}\|\jF{\sigma_h\cdot n_F}\|_{F},
\end{equation*}
where $h_T^{min}=\min\{|F^i|; \  F\in \partial T\cap \F^\Gamma, \  i =1, 2 \}$.

The global error estimator is given by $\eta+\eta_{\Gamma}$, where: 
\begin{equation*}
  \eta =\bigg (\sum_{T\in \mathcal{T}_h}\eta_{T}^2\bigg)^{1/2},
        \qquad \eta_{\Gamma}=\bigg(
         \sum_{F\in \F^\Gamma}\eta_F^2+
         \sum_{T\in\T^\Gamma}\tilde{\eta}_T^2\bigg)^{1/2}.
\end{equation*}
The data approximation terms \Magenta{on a triangle $T\in \mathcal T_h$ and on $\Omega$ are respectively} given by
\begin{equation*}
\Magenta{\epsilon(T)= \frac{h_T}{\sqrt{\delta_T}}\|f-f_h\|_{T}, \quad \epsilon(\Omega)= \bigg(\sum_{T\in \T}\epsilon(T)^2 \bigg)^{1/2}}\quad \text{where}\quad \delta_T= \left\{
\begin{array}{ll}
k_i & \text{if}\ T\in \T^i\backslash\T^\Gamma\\
k_\Gamma &\text {if} \ T\in \T^\Gamma
\end{array}\right..
\end{equation*}
\Magenta{We similarly define $\epsilon(\omega)$ on any subset $\omega\subset\Omega$ that is a union of cells.}
In the following, we establish the reliability and local efficiency of $\eta+\eta_{\Gamma}$.

\section{Reliability}\label{sec:reliab}
Thanks to the fact that $\sigma_h$  strongly satisfies the transmission condition, we are able to establish the reliability of the error estimator $\eta+ \eta_\Gamma$.
As expected\Red{,} when using equilibrated fluxes in a posteriori error estimation, we show that the constant in front of the main part $\eta=\Vert \tau_h\Vert_{\Omega}$ of the global estimator is equal to $1$.

The proof is done in three main steps.

\subsection{Sharp upper bound of the error}\label{subsec:reliab_1}
We begin by establishing a sharp bound for the energy norm of the error. However, at this stage, the upper bound is not yet the global estimator $\eta+\eta_{\Gamma}$.

\begin{lmm}\label{thrm: first_bound_rb}
%Let $u$ and $u_h$ be the solutions of \eqref{eq: continuous_problem} and \eqref{eq: Cutfem_Formulation}, respectively. 
There exists a constant $C>0$ independent of the discretization, the diffusion coefficients, and the interface such that  
\begin{equation*}
|u-u_h|_{1,K,h}\leq \eta +  \inf_{v\in H_0^1(\Omega)}|v-u_h|_{1,K,h} +C\bigg(\sum_{F\in\F^\Gamma}  \eta_F^2+\epsilon (\Omega)^2\bigg)^{1/2}.
\end{equation*}
\end{lmm}
\begin{proof}
We consider $\varphi \in H^1_0(\Omega)$, the unique solution of the following weak problem:  
\begin{equation}\label{eq: aposteriori_phi}
\int_{\Omega}K\nabla\varphi\cdot \nabla v \dx = \int_{\Omega}K\nabla_h u_h\cdot \nabla v \dx \qquad \forall v \in H^1_0(\Omega).
\end{equation}

The triangle inequality gives: 
\begin{equation}\label{eq: triangl_inq_2}
|u-u_h|_{1,K,h}\leq  |u-\varphi|_{1,K,h}+ |\varphi-u_h|_{1,K,h}.
\end{equation}
The auxiliary problem (\ref{eq: aposteriori_phi}) immediately yields that:
\begin{equation}\label{eq:estim_terme_aux}
|\varphi-u_h|_{1,K,h}= \inf_{v\in H_0^1(\Omega)}|v-u_h|_{1,K,h}.%\leq |I_hu_h-u_h|_{1,K,h}\Red{.}
\end{equation}
	
With $\sigma=K\nabla u$ \Magenta{and $\tau_h=K^{-1/2}(\sigma_h-K\nabla_h u_h) $}, one has that:
\begin{equation*}\label{eq: aposteriori_step1}
	\begin{split}
		|u-\varphi|_{1,K,h}^2=&\int_{\Omega}K^{1/2}\nabla(u-\varphi)\cdot(K^{-1/2}\sigma-K^{-1/2} \sigma_h+ \tau_h +K^{1/2}(\nabla_h u_h-\nabla\varphi))\dx\\
		=&\int_{\Omega}\nabla(u-\varphi)\cdot(\sigma-\sigma_h)\dx +\int_{\Omega} K^{1/2}\nabla(u-\varphi)\cdot\tau_h\dx 
		+\int_{\Omega}K\nabla(u-\varphi)\cdot(\nabla_h u_h-\nabla \varphi)\dx.
	\end{split}
\end{equation*}
The last term in the above equality  vanishes thanks to  \eqref{eq: aposteriori_phi}. \Blue{Recalling that $\eta=\Vert \tau_h\Vert_{\Omega}$}, the Cauchy-Schwarz inequality then yields:   
\begin{equation}\label{eq: aposteriori_step2_2}
|u-\varphi|_{1,K,h}^2\leq \left |\int_{\Omega}\nabla(u-\varphi)\cdot(\sigma-\sigma_h)\dx  \right|+\eta|u-\varphi|_{1,K,h}.
\end{equation}
Integrating by parts and using that $(u-\varphi)\in H^1_0(\Omega)$, $\sigma \in H(\div,\Omega)$, $[\![\sigIRT\cdot n_F]\!]=0$ on any non-cut interior edge $F$ and $\div\sigIRT=-f_h $, we further obtain:
\begin{equation}\label{eq: aposteriori_step2_3}
\int_{\Omega}\nabla(u-\varphi)\cdot(\sigma-\sigma_h)\dx =\int_{\Omega}(f-f_h)(u-\varphi)\dx-\sum_{F\in\F^\Gamma} \int_{F}[\![\sigIRT\cdot n_F]\!](u-\varphi)\ds.
\end{equation}

Using the Cauchy-Schwarz inequality, $\nrmL[T]{u-\varphi-\pi_T^0(u-\varphi)}\lesssim h_T \nrmL[T]{\nabla (u-\varphi)}$ for any $T\in \T$ and $k_\Gamma \leq k_i$ for $i=1,2$, we get in a standard way:
\begin{equation}\label{eq: aposteriori_bound_of_f-fh}
\begin{split}
\bigg|	\int_{\Omega}(f-f_h)(u-\varphi)\dx	\bigg|=&\bigg|\sum_{T\in\T}\int_{T}(f-f_h)(u-\varphi -\pi_T^0(u-\varphi))\dx	\bigg|\\
\lesssim  &\sum_{i=1}^2\sum_{T\in\T^i \backslash \T^\Gamma}\frac{h_T}{\sqrt{k_i}}\|f-f_h\|_{T}\|k_i^{1/2}\nabla(u-\varphi)\|_{T}\\
&+\sum_{T\in \T^\Gamma}\frac{h_T}{\sqrt{k_{\Gamma}}}\|f-f_h\|_{T}\|k_\Gamma^{1/2}\nabla(u-\varphi)\|_{T}\\
\lesssim &|u-\varphi|_{1,K,h}\, \epsilon(\Omega).
\end{split}
\end{equation}

Let $F\in \F^\Gamma$. By definition of the space $\IRT^0$, we have $\pi_F^0\jF{\sigIRT\cdot n_F}=0$, hence:
\begin{equation*}
\bigg| \sum_{F\in\F^\Gamma} \int_{F}\jF{\sigIRT\cdot n_F}(u-\varphi)\ds  \bigg|
\leq \sum_{F\in\F^\Gamma} \frac{\sqrt{h_F}}{\sqrt{k_\Gamma}}\|\jF{\sigIRT\cdot n_F}\|_{F}\frac{\sqrt{k_\Gamma}}{\sqrt{h_F}}\|u-\varphi-\pi_F^0(u-\varphi)\|_{F}. %-\pi^0_F\jF{\sigIRT\cdot n_F}.
\end{equation*}
We have that
$h_F^{-1/2}\|u-\varphi-\pi_{\Blue{F}}^0(u-\varphi)\|_{F}\lesssim  |u-\varphi|_{1,T_F^+\cup T_F^-}$. \Blue{The proof of this bound is standard: we first use the minimal distance property of $\pi_F^0$, which yields $\|u-\varphi-\pi_{F}^0(u-\varphi)\|_{F}\le \|u-\varphi-\pi_{T}^0(u-\varphi)\|_{F}$ for any triangle $T$ such that $F\subset \partial T$, then the well-known trace inequality $h_F^{-1/2}\|v\|_{F}\lesssim h_T^{-1}\|v\|_{T}+|v|_{1,T}$ for any $v\in H^1(T)$ and finally, the interpolation error estimate $\|u-\varphi-\pi_{T}^0 (u-\varphi)\|_{T}\lesssim  h_T|u-\varphi|_{1,T}$. Next,} using the Cauchy-Schwarz inequality and $k_\Gamma\leq k_i$ for $i=1,2$, we obtain:
\begin{equation}\label{eq: bound _weak continuity}
\bigg| \sum_{F\in\F^\Gamma} \int_{F}\jF{\sigIRT\cdot n_F}(u-\varphi)\ds\bigg|
\lesssim |u-\varphi|_{1,K,h} 
\bigg(\sum_{F\in\F^\Gamma}  \eta_F^2\bigg)^{1/2}.
\end{equation}

Gathering together \eqref{eq: aposteriori_step2_2}, \eqref{eq: aposteriori_step2_3}, \eqref{eq: aposteriori_bound_of_f-fh} and \eqref{eq: bound _weak continuity}, we obtain the following bound:
\begin{equation*}\label{eq: aposteriori_bound_1}
|u-\varphi|_{1,K,h}\leq \eta + C\bigg(\sum_{F\in\F^\Gamma}  \eta_F^2+\epsilon (\Omega)^2\bigg)^{1/2},
\end{equation*}
which together with (\ref{eq: triangl_inq_2}) and (\ref{eq:estim_terme_aux}) ends the proof.
\end{proof}

Note that 
\begin{equation}\label{eq: inf_bound}
\displaystyle \inf_{v\in H_0^1(\Omega)}|v-u_h|_{1,K,h}\leq |I_hu_h-u_h|_{1,K,h},	
\end{equation}
for any interpolation $I_hu_h$ of $u_h$ belonging to $H^1_0(\Omega)$. Therefore, the next step consists \Red{of} constructing $I_hu_h\in H^1_0(\Omega)$ such that the interpolation error is bounded by
$\eta_\Gamma$. 

\subsection{A continuous approximation of the CutFEM solution}\label{subsec:interp}
We focus on the cut triangles since on $T\in\mathcal T_h\backslash \mathcal T_h^\Gamma$, we simply take $I_hu_h=u_h$.

Let $T\in \mathcal T_h^\Gamma$. We denote the triangle $T$ by $\triangle A_1A_2A_3$ and let $M\in A_1A_2$,  $N\in A_1A_3$ be the intersection points with $\Gamma$. \Blue{The case where $T$ is cut at one of its vertices into two sub-triangles can be treated in a similar manner and leads to analogous results (see Remark \ref{rk:1}). Therefore, in what follows, we only detail the case where one of the cut parts is quadrilateral, which is more challenging.} 

\begin{figure}[ht]
	\centering
	\includegraphics[width=0.28\textwidth]{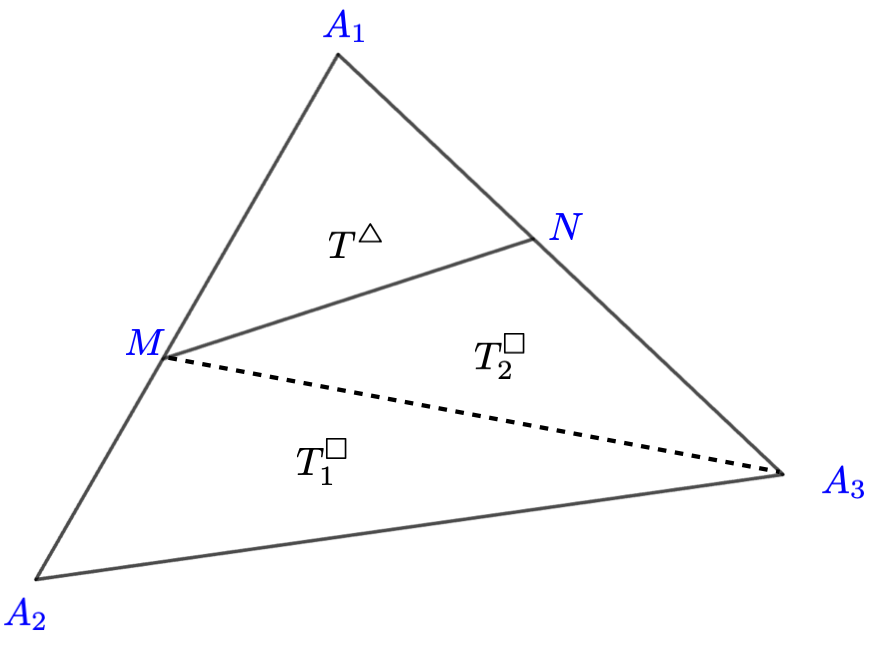}
	\caption{Notation on a cut element in view of the interpolation}
	\label{Fig: cut_I_h_TRiangle}
\end{figure}

For simplicity of notation, let $T^\triangle$ and $T^\square$ denote the triangular and quadrilateral parts of $T$, respectively, and assume that $T^\triangle=T^1$, $T^\square= T^2$. Since the ratios   $|A_2M| /|A_3N|$ and  $|A_3N| / |A_2M|$ cannot  blow up simultaneously, we assume, without loss of generality, that $|A_2M| / |A_3N|$ is bounded and then cut $T^\square$ into two triangles using $A_3M$, cf. Figure \ref{Fig: cut_I_h_TRiangle} (otherwise, use $A_2N$).

We define the interpolation operator $I_h$ as follows: $I_hu_h$ is a linear function on each sub-triangle $\triangle A_1MN=:T^\triangle$, $\triangle MA_2A_3=:T^\square_1$ and $\triangle MNA_3=:T^\square_2$, and
\begin{equation}\label{eq: def_I_h_0}
I_hu_h(A_1)=u_{h,1}(A_1),\quad I_hu_h(A_2)=u_{h,2}(A_2), \quad I_hu_h(A_3)=u_{h,2}(A_3),
\end{equation}
\begin{equation}\label{eq: def_I_h}
I_hu_h(M)=\{u_h\}^*(M),\quad I_hu_h(N)=\{u_h\}^*(N).
\end{equation}

This choice guarantees that $I_h u_h\in H^1(T)$, since it is continuous at the intersection points $M$, $N$ and $A_3$ of the sub-triangles. Since $[\![{I_hu_h}]\!]_F=0$ for any edge $F\in \mathcal F_h$, we have $I_hu_h\in H^1_0(\Omega)$.

%\begin{remark}
%A naive way to define $I_hu_h$ is to take $I_hu_h\in P^1(T)$ that satisfies (\ref{eq: def_I_h_0}). In this case, we can establish the following estimate: 
%\begin{equation*}
%|I_h u_h-u_h|_{1,K,h}^2\lesssim  \sum_{T\in \T^\Gamma} k_{max}\sum_{N\in \mathcal N_T}\jump{u_h}^2(N)\lesssim \sum_{T\in \T^\Gamma}\frac{k_{max}}{h_T^2}\|\jump{u_h}\|_T^2 
%\end{equation*}
%but we couldn't bound it further by $\eta_\Gamma$. 
%\end{remark}

In view of estimating the interpolation error $|I_hu_h-u_h|_{1,K,T}$, we first write it as: 
\begin{equation}\label{eq: decomp}
|I_h u_h-u_h|_{1,K,T}^2=|I_hu_h-u_{h,1}|_{1,K,T^\triangle}^2+\sum_{l=1}^2|I_hu_h-u_{h,2}|_{1,K,T^\square_l}^2,
\end{equation}
\Red{where for any $\tilde T\in \{T^\triangle,T^\square_1,T^\square_2\}$ and any $v\in H^1(\tilde T)$, we set $|v|_{1,K,\tilde T}=\|K^{1/2}\nabla v\|_{\tilde T} $.} Let now $\tilde T\in \{T^\triangle,T^\square_1,T^\square_2\}$ and assume $\tilde T\Blue{\subset} \Omega^i$, with $i\in\{1,2\}$. Then we have, thanks to \eqref{eq: def_I_h_0} and \eqref{eq: def_I_h}, that:
\begin{equation*}
(I_hu_h-u_h)_{|\tilde T}= (-1)^i\omega_i[u_h](M) \varphi_M +(-1)^i\omega_i[u_h](N)\varphi_N, 
\end{equation*}
where $\varphi_M, \varphi_N\in P^1(\tilde T)$ are the nodal basis functions associated to $M$ and $N$. We express them below in terms of the barycentric coordinates $(\lambda_j)_{1\le j\le 3}$ of $T$. 
\begin{lmm}
One has that:
\begin{equation}\label{eq: phi_M_T_delta}
(\varphi_M)_{|T^\triangle}=\frac{|F_3|}{|A_1M|}\lambda_2, \quad (\varphi_N)_{|T^\triangle}=\frac{|F_2|}{|A_1N|}\lambda_3,
\end{equation} 
\begin{equation}\label{eq: phi_M_T_squar_1}
(\varphi_M)_{|T^\square_1}=\frac{|F_3|}{|A_2M|}\lambda_1,	\quad (\varphi_N)_{|T^\square_1}=0,
\end{equation}
\begin{equation}\label{eq: phi_M_T_squar_2}
(\varphi_M)_{|T^\square_2}=\frac{|F_3|}{|A_1M|}\lambda_2,	 \quad (\varphi_N)_{|T^\square_2}= \frac{|F_2|}{|A_3N|} \lambda_1 - \frac{|F_2|}{|A_3N|} \frac{|A_2M|}{|A_1M|}\lambda_2,
\end{equation}
where $F_j$ denotes the edge opposite to the node $A_j$, for $1\leq j\leq 3$.
\end{lmm}
\begin{proof}
On $T^\triangle$, we have $\varphi_M(A_1)=\varphi_M(N)=0$ and ${\lambda_2}_{|F_2}=0$, so $\varphi_M= c\lambda_2$ with $c\in\mathbb R$. From the relation $\lambda_2(M)= |A_1M|/|F_3|$, we obtain the desired expression of $ (\varphi_M)_{|T^\triangle} $. We similarly get the expression of $(\varphi_N)_{|T^\triangle}$, so (\ref{eq: phi_M_T_delta}) is established.

On $T^\square_1$, we immediately get $(\varphi_N)_{|T^\square_1}=0$. The expression of $\varphi_M$ is obtained in a similar way to the previous case, thus yielding \eqref{eq: phi_M_T_squar_1}. Finally, on $T^\square_2$, we similarly get $\varphi_M$ while $\varphi_N$ is obtained by decomposing it in the barycentric basis: 
\begin{equation*}
(\varphi_N)_{|T^\square_2}=a \lambda_1+b\lambda_2+c\lambda_3,\qquad a,b,c\in\mathbb R.
\end{equation*}
From $\varphi_N(A_3)=0$, we get $c=0$. Moreover, we have $\varphi_N(M)=a\lambda_1(M)+b\lambda_2(M)=0 $ and $\varphi_N(N)= a \lambda_1(N)=1$. Using next 
\[
\lambda_1(M)=\frac{|A_2M|}{|F_3|}, \quad \lambda_1(N)=\frac{|A_3N|}{|F_2|}, \quad \lambda_2(M)=\frac{|A_1M|}{|F_3|},
\]
we compute $a$ and $b$ and obtain 
%\[ a=\frac{1}{\lambda_1(N)}=\frac{|F_2|}{|A_3N|},\quad  b= -\frac{|F_2|}{|A_3N|}\frac{\lambda_1(M)}{\lambda_2(M)}=  -\frac{|F_2|}{|A_3N|} \frac{|A_2M|}{|A_1M|},
%\]
the relation \eqref{eq: phi_M_T_squar_2}.  
\end{proof}

\subsection{Estimate of the interpolation error}\label{subsec:estim_interp}
We can now bound the interpolation error $|I_hu_h-u_h|_{1,K,T}$ on any cut triangle $T$. We begin by establishing an auxiliary estimate. For this purpose, let $\alpha :=[u_h](M)$ and $\beta:=[u_h](N)$, as well as
\begin{equation}\label{eq:def_expr}
E_\alpha:= \frac{|T^\triangle|}{|A_1M|^2}+ \frac{|T^\square_1|}{|A_2M|^2}+\frac{|T^\square_2|}{|A_1M|^2}, \qquad E_\beta:=\frac{|T^\triangle|}{|A_1N|^2}+ \frac{|T^\square_2|}{|A_3N|^2}\frac{|F_3|^2}{|A_1M|^2}.
\end{equation}
\begin{lmm}\label{lem:expressions_0}
One has that:
    \begin{equation*}\label{eq: bound interpolation E_A_E_b}
		|I_hu_h-u_h|_{1,K,T}^2\lesssim k_\Gamma (\alpha^2 E_\alpha+\beta^2 E_\beta).
    \end{equation*}
\end{lmm}
\begin{proof}
We use \eqref{eq: decomp} and bound the error on each sub-triangle $\tilde T$ of $T$. 

First, let $\tilde  T=T^\triangle$. Using that $k_\Gamma=k_i\omega_i$ and $0< \omega_i< 1$ for $i=1, 2$, we get
\begin{equation*}
		|I_hu_h-u_{h,1}|_{1,K,T^\triangle}^2 \leq  k_\Gamma |W_h^\triangle|_{1,T^\triangle}^2
\end{equation*}
where $W_h^\triangle$ belongs to $P^1(T^\triangle)$ and satisfies $W_h^\triangle(M)= -\alpha$, $W_h^\triangle(N)=-\beta$ and $W_h^\triangle(A_1)=0$. Hence, $W_h^\triangle=-\alpha \varphi_M-\beta \varphi_N$. From \eqref{eq: phi_M_T_delta}, combined with the relations
%we get:
%\begin{equation*}
%|W_h^\triangle|_{1,T^\triangle}^2\lesssim \alpha^2\frac{|F_3|^2}{|A_1M|^2}| \lambda_2|_{1,T^\triangle}^2+\beta^2\frac{|F_2|^2}{|A_1N|^2}|\lambda_3|_{1,T^\triangle}^2. 
%\end{equation*}
\begin{equation}\label{eq:barycentric_ppty}
|\lambda_j|_{1,T^\triangle}^2=\frac{|T^\triangle|}{|T|}|\lambda_j|_{1,T}^2,\qquad  |\lambda_j|_{1,T}\lesssim 1\quad (1\leq j\leq 3),\qquad |T|\approx h_T^2
\end{equation}
and $|F_j|\le h_T \,(1\le j\le3)$, we obtain the following bound:  
\begin{equation}\label{eq: W_h_0}
|I_hu_h-u_{h,1}|_{1,K,T^\triangle}^2\lesssim k_{\Gamma}\bigg(\alpha^2\frac{|T^\triangle|}{|A_1M|^2}+\beta^2\frac{|T^\triangle|}{|A_1N|^2}\bigg). 
\end{equation}

Next, let $\tilde T=T^\square_1$. Similarly, one has $|I_hu_h-u_{h,2}|_{1,K,T^\square_1}^2 \leq  k_\Gamma |W_h^1|_{1,T^\square_1}^2$ where $W_h^1\in P^1(T^\square_1)$ and now satisfies: $W_h^1(M)= \alpha$, $W_h^1(A_2)= W_h^1(A_3)=0$. Hence, $W_h^1=\alpha \varphi_M$ and from \eqref{eq: phi_M_T_squar_1} and \eqref{eq:barycentric_ppty} we obtain:
\begin{equation}\label{eq: W_h_1}
|I_hu_h-u_{h,2}|_{1,K,T^\square_1}^2  %k_\Gamma\alpha^2 \frac{|T^\square_1|}{|T|} \frac{|F_3|^2}{|A_2M|^2} |\lambda_1|_{1,T}^2
\lesssim  k_\Gamma \alpha^2\frac{|T^\square_1|}{|A_2M|^2}.
\end{equation}

Finally, for  $\tilde T=T^\square_2$ we similarly have that
\begin{equation*}
|I_hu_h-u_{h,2}|_{1,K,T^\square_2}^2 \leq  k_\Gamma |W_h^2|_{1,T^\square_2}^2 
\end{equation*}
where $W_h^2\in P^1(T^\square_2)$ with $W_h^2(M)= \alpha$, $ W_h^1(N)=\beta$ and $ W_h^2(A_3)=0$, which yields $W_h^2= \alpha \varphi_M+ \beta \varphi_N$. From \eqref{eq: phi_M_T_squar_2}, we get  
\begin{equation*}
|W_h^2|_{1,T^\square_2}^2\lesssim \alpha^2 \frac{|F_3|^2}{|A_1M|^2}\frac{|T^\square_2|}{|T|} |\lambda_2|_{1,T}^2+ \beta^2 \frac{|F_2|^2}{|A_3N|^2}\frac{|T^\square_2|}{|T|}\bigg (|\lambda_1|_{1,T}^2+\frac{|A_2M|^2}{|A_1M|^2}|\lambda_2|_{1,T}^2
\bigg).
\end{equation*}
Using again (\ref{eq:barycentric_ppty}), we obtain the following inequality:  
\begin{equation*} 	
\frac{|F_2|^2}{|A_3N|^2}\frac{|T^\square_2|}{|T|} \bigg(|\lambda_1|_{1,T}^2+\frac{|A_2M|^2}{|A_1M|^2} |\lambda_2|_{1,T}^2\bigg) \lesssim \frac{|T^\square_2|}{|A_3N|^2}\bigg(1+\frac{|A_2M|^2}{|A_1M|^2}\bigg)
\lesssim \frac{|T^\square_2|}{|A_3N|^2}\frac{|F_3|^2}{|A_1M|^2}.
\end{equation*}
Hence, we finally deduce that 
\begin{equation}\label{eq: W_h_2}
|I_hu_h-u_{h,2}|_{1,K,T^\square_2}^2\lesssim k_{\Gamma}\bigg(\alpha^2 \frac{|T^\square_2|}{|A_1M|^2} + \beta^2 \frac{|T^\square_2|}{|A_3N|^2}\frac{|F_3|^2}{|A_1M|^2}\bigg).
\end{equation} 
Gathering together the estimates \eqref{eq: W_h_0}, \eqref{eq: W_h_1} and \eqref{eq: W_h_2} in \eqref{eq: decomp} ends the proof.
\end{proof}
	
The next step consists of bounding the terms $E_\alpha$ and $E_\beta$ in the best possible way with respect to the geometry of the cut cell $T$. A first estimate is given in the following lemma, based on equivalent expressions of these terms.
\begin{lmm}\label{lem:expressions}
One has that:
\begin{equation*}\label{eq: E_alpha_beta}
E_\alpha=\frac{|T|}{|A_1M||A_2M|},\qquad E_\beta\le \frac{|T|}{|A_1N||A_3N|}+\frac{|T|}{|F_2||A_1M|}\frac{|A_2M|}{|A_3N|}.
\end{equation*}
\end{lmm}
\begin{proof}
We start by developing the expressions of $E_\alpha$ and  $E_\beta$ introduced in (\ref{eq:def_expr}). For this purpose, it is useful to note that
\begin{equation*}
\frac{|T^\triangle|}{|T|}%=\frac{|A_1M||A_1N|\sin \hat A_1 }{|F_3||F_2|\sin \hat A_1}
= \frac{|A_1M||A_1N| }{|F_3||F_2|},\quad 
\frac{|T^\square_1|}{|T|}%=\frac{|A_2M||F_1|\sin \hat A_2 }{|F_3||F_1|\sin \hat A_2}
= \frac{|A_2M|}{|F_3|}, \quad \frac{|T^\square_2|}{|T|}=1-
\frac{|T^\square_1|}{|T|} -\frac{|T^\triangle|}{|T|}.
\end{equation*}
%\begin{equation*}
%\frac{|T^\square_2|}{|T|}= \frac{|T|-|T^\square_1|-|T^\triangle|}{|T|}=1-\frac{|A_2M|}{|F_3|}-  \frac{|A_1M||A_1N| }{|F_3||F_2|}.
%\end{equation*}
Using $|A_1N|+|A_3N|=|F_2|$ and $|A_1M|+|A_2M|=|F_3|$, we further obtain
\begin{equation*}
\frac{|T^\square_2|}{|T|}=\frac{|A_1M|}{|F_3|}\bigg (1- \frac{|A_1N| }{|F_2|}\bigg )= \frac{|A_1M|}{|F_3|}\frac{|A_3N| }{|F_2|}.
\end{equation*}

Substituting $|T^\triangle|$, $|T^\square_1|$ and $|T^\square_2|$ from the previous relations in (\ref{eq:def_expr}), we get after some straightforward computation:
\begin{equation*}
\begin{split}
E_\alpha=&|T|\bigg(\frac{1}{|A_1M||F_3|}+\frac{1}{|A_2M||F_3|}\bigg)=\frac{|T|}{|A_1M||MA_2|},\\
E_\beta=&\frac{|T|}{|F_2|}\bigg(\frac{|A_1M|}{|A_1N||F_3|}+\frac{|F_3|}{|A_1M||A_3N|}\bigg).
\end{split}
\end{equation*} 
In the expression of $E_\beta$, we next use $|A_1M|\leq |F_3|$ and substitute $|F_3|=|A_1M|+|A_2M|$ to obtain: 
\begin{equation*}
E_\beta \leq \frac{|T|}{|F_2|}\bigg(\frac{1}{|A_1N|}+\frac{1}{|A_3N|}+\frac{|A_2M|}{|A_1M||A_3N|}\bigg) = \frac{|T|}{|A_1N||A_3N|}+\frac{|T|}{|F_2||A_1M|}\frac{|A_2M|}{|A_3N|}.
\end{equation*}
This ends the lemma's proof.
\end{proof}

Next, we bound $E_\alpha$ and $E_\beta$ by a simpler expression. 
\begin{lmm}\label{lem:bound_E_final}
One has that:
\begin{equation*}\label{eq: E_final}
E_\alpha\le\frac{h_T}{h_T^{min}}, \quad E_\beta \lesssim  \frac{h_T}{h_T^{min}}.
\end{equation*} 
\end{lmm}
\begin{proof}
Since $|A_1M|+|A_2M|=|F_3|$, we have that:  
$$\displaystyle |A_1M||A_2M|\geq \min\{|A_1M|,|A_2M|\}\frac{|F_3|}{2}.$$
So from Lemma \ref{lem:expressions} we immediately deduce the bound for $E_\alpha$:  
\begin{equation*}
E_{\alpha}\leq \frac{2|T|}{h_T^{min}|F_3|}\le \frac{h_T}{h_T^{min}},
\end{equation*}
recalling that $h_T^{min}= \min\{|A_1M|,|A_2M|, |A_1N|,|A_3N|\}$. Using now $|A_1N|+|A_3N|=|F_2|$, we similarly get
\begin{equation*}
 E_\beta \le  \frac{h_T}{h_T^{min}}+ \frac{|T|}{|F_2|h_T^{min}} \frac{|A_2M|}{|A_3N|}\le \frac{h_T}{h_T^{min}}+ \frac{2h_T}{h_T^{min}} \frac{|A_2M|}{|A_3N|}.
\end{equation*}
Our initial choice of cutting $T^\square$ into two sub-triangles ensures that $|A_2M|/|A_3N|$ is bounded (otherwise, it suffices to change the roles between $ M $ and $N$). 
%In fact, the only case where $E_\beta$ could blow up occurs when $M$ is close to $A_1$ and $N$ close to $A_3$. %(see Figure \ref{Fig: cut_cases_discuss}-(d)). However, in this case, we can change the roles between $N$ and $M$, that is cut $T^\square$ by the node $N $ instead of $M$. 
So we also get $E_\beta \lesssim  h_T/h_T^{min}$, which concludes the proof.
\end{proof}

\Blue{
\begin{remark}\label{rk:1}
Lemma \ref{lem:bound_E_final} remains unchanged in the case where $T$ is cut by $\Gamma$ at one of its vertices - let's say $A_3=N$ - into two sub-triangles. One can then choose either of the sub-triangles as $T^{\Delta}$, let's say $T^{\Delta}=A_1A_3M$. Then we have that $\varphi_N=\lambda_3$, $(\varphi_M)_{|T^{\Delta}}=\dfrac{|F_3|}{|A_1M|}\lambda_2$ and  $(\varphi_M)_{|T^{\square}}=\dfrac{|F_3|}{|A_2M|}\lambda_1$.
The proofs of Lemmas \ref{lem:expressions_0}, \ref{eq: E_alpha_beta} and \ref{eq: E_final} follow the same idea as in the quadrilateral case, though they are simpler. We thus get 
$$E_{\alpha}=1,\qquad E_{\beta}=\dfrac{|T^\triangle|}{|A_1M|^2}+\dfrac{|T^\square|}{|A_2M|^2}\le \dfrac{h_T}{2}\bigg(\dfrac{1}{|A_1M|}+\dfrac{1}{|A_2M|}\bigg)=\dfrac{h_T|F_3|}{2|A_1M||A_2M|}\le \dfrac{h_T}{h_T^{min}}.$$ 
\end{remark}
}

We can now establish the interpolation error estimate in a cut cell. 
\begin{lmm}\label{thm:final_interp_est}
Let $T\in \T^\Gamma$. Then one has  $|I_hu_h-u_h|_{1,K,T}\lesssim \tilde{\eta}_T$.
\end{lmm}
\begin{proof}
Lemmas \Magenta{\ref{lem:expressions_0}} and \ref{lem:bound_E_final} yield:
\begin{equation*}
|I_hu_h-u_h|_{1,K,T}^2\lesssim \frac{h_Tk_\Gamma}{h_T^{min}} (\alpha^2+\beta^2)\lesssim \frac{h_Tk_\Gamma}{h_T^{min}|\Gamma_T|}\|[u_h]\|_{\Gamma_T}^2=\tilde{\eta}_T^2
\end{equation*}
which is the desired result.
\end{proof}

By summing the estimate of Lemma \ref{thm:final_interp_est} upon the cut cells, we immediately obtain from (\ref{eq: inf_bound}) the next bound.
 
\begin{lmm}\label{lem: bound_inf}%For any $T\in \T^\Gamma$
One has that $\displaystyle \inf_{v\in H_0^1(\Omega)}|v-u_h|_{1,K,h}\lesssim  \eta_\Gamma$.	\end{lmm} 
 
Finally, thanks to Lemma \ref{thrm: first_bound_rb} and Lemma \ref{lem: bound_inf}, we deduce the reliability bound.

\begin{thrm}[Reliability]
There exists a constant $C>0$ independent of the mesh, the coefficients and the interface such that   
\begin{equation}\label{eq:reliab}
	|u-u_h|_{1,K,h}\leq \eta +C\left( \eta_{\Gamma}+\epsilon(\Omega)\right).
\end{equation}
\end{thrm}

This theorem shows that the weighted $H^1$ semi-norm of the error is bounded by \Blue{the global estimator $\eta+\eta_{\Gamma}$ plus the higher order term $\epsilon(\Omega)$, with a reliability constant which is independent of the mesh size, the diffusion coefficients, and the mesh/interface configuration. In addition, the constant in front of the main estimator $\eta$ is equal to $1$, as expected.}

\begin{remark}
\Blue{The price to pay for obtaining a reliability constant independent of the mesh/interface geometry was to weight the estimator $\tilde {\eta}_T$ on a cut triangle by $\sqrt{h_T}/\sqrt{h_T^{min}}$, which arises from our definition of $I_h$. Although this factor may become unbounded, it does not seem to influence the estimator's behavior in the numerical tests. If one is able to construct a continuous interpolation $I_hu_h$ of $u_h$, for which the interpolation error is fully robust  with respect to the mesh/interface configuration, then the factor $\sqrt{h_T}/\sqrt{h_T^{min}}$ would disappear.}
\end{remark}

\section{Properties of a piecewise Raviart-Thomas space on the cut cells}\label{sec:decomp_RT}

To establish the local efficiency of $\eta_T=\Vert K^{-1/2}(\sigma_h-K\nabla_hu_h)\Vert_{T}$ on a cut triangle $T\Red{ \in\mathcal T_h^{\Gamma}}$, we first need some properties of the finite element space to which $\sigma_h-K\nabla_hu_h$ belongs. Note that although $\sigma_h$ belongs to the immersed Raviart-Thomas space $\IRT^0(T)$,  $\sigIRT -K\nabla_h u_h$ is only piecewise $\RT^0$ on a cut \Red{triangle $T$, i.e. its restriction to each part $T^i=T\cap \Omega^i$  belongs to $\RT^0$; this is} due to the discontinuity of the diffusion coefficients $ K $ and of the approximate solution $ u_h$.

In this section, we provide some properties of this space in relation \Red{to} the immersed one, in particular\Red{,} we construct and bound specific shape functions which will be used in Section \ref{sec:eff} to decompose $\sigma_h-K\nabla_hu_h$.

\subsection{Auxiliary results}\label{subsec:aux_results}
We recall (cf. \cite{raviart1977mixed}) the shape functions of the standard Raviart-Thomas space $\RT^0(T)$ on a triangle $T$ of nodes $(A_j)_{1\leq j\leq 3}$. The basis function $\Lambda_{T,j}$, associated to the edge $F_j$ opposite to $A_j$, satisfies:
\begin{equation}\label{eq: bound base RT}
\Lambda_{T,j}(x)=\frac{|F_j|}{2|T|} \overset{\longrightarrow}{A_jx},\quad |\Lambda_{T,j}(x)|\lesssim 1,\quad \nrmL[T]{\Lambda_{T,j}}\lesssim h_T\quad (1\leq j\leq 3).
\end{equation}
The local Raviart-Thomas interpolation operator is denoted by $I_T$.

On a cut triangle $T\in \T^\Gamma$, we agree to denote by $n_\Gamma^\triangle$ the unit normal vector to $\Gamma$ pointing from $T^\triangle$ to $T^\square$, and by $t_\Gamma^\triangle$ the unit tangent vector oriented by a $90^\circ$ clockwise rotation of $n_\Gamma^\triangle$. We consider the following functions $\omega$ and $\rho$ on $T$:  
\begin{equation}\label{eq:omega_rho}
\omega =\left\{ \begin{array}{ll}
t_\Gamma^\triangle &  \text{ in} \ T^\triangle\\
			0 & \text{ in} \  T^\square
		\end{array}\right.,\quad 
\rho=\left\{ \begin{array}{ll}
		n_\Gamma^\triangle & \text{ in} \ T^{\triangle}\\
		0 &  \text{ in} \ T^{\square}.\\
	\end{array}\right.
\end{equation}

The function $\omega$ was already introduced in \cite{IRT}, where it was used to prove the unisolvency of $\IRT^0(T)$. The author established, under the hypothesis that \Blue{the cut triangle} $T$ is acute, the following bound:
\begin{equation}\label{eq: bound_omega_2}
0\leq( I_T \omega \cdot t_\Gamma^\triangle)(x)\leq 1 ,\quad \forall x\in T.
\end{equation}
%The condition $\alpha_{max}\leq \pi/2$, used in \cite{IRT} in the proof of the unisolvance of $\IRT^0$, is not restrictive since it can be achieved by eventually refining the triangle. 

We can easily prove a similar estimate for the new function $\rho$, which is:
\begin{equation}\label{eq: bound_rho}
|(I_T\rho\cdot t_\Gamma^\triangle)(x)|\lesssim 1\quad  \forall x\in T.
\end{equation}
Indeed, thanks to $|\Lambda_{T,j}(x)|\lesssim 1 $ and $ |\rho(x)|\leq 1 $ for all $ x\in T $, we get
\begin{equation*}
|(I_T\rho) (x)\cdot t_\Gamma^\triangle|\leq  \sum_{j=1}^3 |N_{T,j}(\rho)|| \Lambda_{T,j}(x)|\lesssim 	 \sum_{j=1}^3  \frac{1}{|F_j|}\int_{F_j} |\rho\cdot n_T|ds \lesssim 1.
\end{equation*}
%\begin{figure}[H]
%	\centering
%	\begin{subfigure}[ht]{0.45\textwidth}
%		\includegraphics[width=0.8\textwidth]{Figures/case1.png}
%		\caption{Case 1 ($T^1=T^\triangle$)}
%	\end{subfigure}
%\begin{subfigure}[ht]{0.45\textwidth}
%	\includegraphics[width=0.8\textwidth]{Figures/case2.png}
%	\caption{Case 2 ($T^1=T^\square$)}
%\end{subfigure}
%	\caption{Geometry of a cut triangle }
%	\label{fig:cases cut}
%\end{figure}  

Next, we consider the following functions on the cut triangle $T$:
\begin{equation}\label{eq:phi_t_n}
\phi_t=\omega-I_T\omega, \qquad \phi_n=\rho- I_T\rho
\end{equation}
which are clearly piecewise $\RT^0$ and satisfy: $N_{T,j}(\phi_t)=N_{T,j}(\phi_n)=0$ for $1\le j\le 3$. They also satisfy the following properties.

\begin{lmm}\label{lem:bound phi_t_n}
Let $F_j^\triangle=F_j\cap \partial T^\triangle$ for $j=2,3$. One has that:
\begin{equation*}
\nrmL[T]{\phi_t}\lesssim   \sum_{j=2}^3|F_j^\triangle|, \qquad \nrmL[T]{\phi_n}\lesssim   \sum_{j=2}^3|F_j^\triangle|.
\end{equation*}
\end{lmm}
\begin{proof}
From the definitions (\ref{eq:omega_rho}) of $\omega$ and $\rho$, we immediately get that 
\begin{equation}\label{eq: bound_eta_omega_L2}
\nrmL[T]{\rho}^2= \nrmL[T]{\omega}^2=|T^{\triangle}|\leq\frac{1}{2}|F_2^\triangle||F_3^\triangle|
\leq \frac{1}{4} \bigg(\sum_{j=2}^3|F_j^\triangle|\bigg)^2.
\end{equation}
By definition of the Raviart-Thomas interpolation operator $ I_T $, we get thanks to \eqref{eq: bound base RT}:
\begin{equation*}
\nrmL[T]{I_T\rho}\lesssim  h_T\sum_{j=1}^{3}|N_{T,j}(\rho)|,\qquad 
\nrmL[T]{I_T\omega}\lesssim  h_T\sum_{j=1}^{3}|N_{T,j}(\omega)|,
\end{equation*}
where $ N_{T,1}(\rho)=N_{T,1}(\omega)=0$ and where we have, for $2\leq j\leq 3$, that
\begin{equation*}
\begin{split}
|N_{T,j}(\rho)|\leq \frac{1}{|F_j|}\int_{F_j}|\rho\cdot n_{T}|\ds  = \frac{1}{|F_j|}\int_{F_j^{\vartriangle}}|n_\Gamma^\triangle\cdot n_{T}|\ds \leq \frac{|F_j^{\vartriangle}|}{|F_j|},\\\
|N_{T,j}(\omega)|\leq \frac{1}{|F_j|}\int_{F_j}|\omega\cdot n_{T}|\ds  = \frac{1}{|F_j|}\int_{F_j^{\vartriangle}}|t_\Gamma^\triangle\cdot n_{T}|\ds \leq \frac{|F_j^{\vartriangle}|}{|F_j|}.
\end{split}
\end{equation*}
Hence, we end up with $\nrmL[T]{I_T\omega}\lesssim \displaystyle\sum_{j=2}^3|F_j^\triangle|$ and $\nrmL[T]{I_T\rho}\lesssim \displaystyle\sum_{j=2}^3|F_j^\triangle|$, which together with \eqref{eq:phi_t_n} and \eqref{eq: bound_eta_omega_L2} yield the result.
\end{proof}

\Red{It is useful to recall that $div (\phi_t )$ and $div (\phi_n)$ are constant on each $T^i$ ($1\le i\le 2$), and that $[\phi_t\cdot n_{\Gamma}]$ and $[\phi_n\cdot n_{\Gamma}]$ are constant on $\Gamma_T$, thanks to the properties of the Raviart-Thomas space $\mathcal {RT}^0$; we also recall that $x_{\Gamma}$ is the arbitrary point of $\Gamma_T$ considered in the definition of $\mathcal {IRT}^0(T)$, cf. Section \ref{sec: A posteriori_estim}.
}

\begin{lmm}\label{lem:ppty_phi}
\Red{Let $T\in \mathcal T_h^{\Gamma}$.}One has that:
\begin{eqnarray*}
\div({\phi_t}_{|T^1})=\div({\phi_t}_{|T^2}),\quad [\phi_t\cdot n_\Gamma]=0,\quad [\phi_t\cdot t_\Gamma](x_\Gamma)=1,\\
\div({\phi_n}_{|T^1})=\div({\phi_n}_{|T^2}),\quad [\phi_n\cdot n_\Gamma]=1,\quad [\phi_n\cdot t_\Gamma](x_\Gamma)=0.
\end{eqnarray*}
\end{lmm}
\begin{proof}
The relations concerning $\phi_t$ were already established in \cite{IRT}. The proof for $\phi_n$ is similar. Since $\div(I_T\rho)$ is constant on $T$ and $\rho$ is piecewise constant, we get the first equality. Then, using that $I_T\rho$ is continuous across $\Gamma_T$ and that $n_{\Gamma}^\Delta\cdot t_\Gamma=0$, we obtain the last \Red{equality on $\Gamma_T$, and in particular at the point $x_{\Gamma}$.} Considering then the cases $T^\Delta=T^1$ and $T^\Delta=T^2$, we also get the second relation: $[\phi_n\cdot n_\Gamma]=[\rho\cdot n_\Gamma]=(n_\Gamma^\Delta \cdot n_\Gamma)^2=1$ \Red{on $\Gamma_T$}.     
\end{proof}

Thus, $\phi_t$ and $\phi_n$ belong to the following subspace of $ \mathcal{RT}^0(T^1)\times \mathcal{RT}^0(T^2)$:
\begin{equation*}
\mathcal E(T):= \left\{\tau=(\tau_1,\tau_2) \in   \mathcal{RT}^0(T^1)\times \mathcal{RT}^0(T^2);\quad  \div\tau_1 =\div\tau_2  \right\}.
\end{equation*}
To conclude this subsection, let us introduce the linear operator $\mathcal R_T :\mathcal E (T)\longrightarrow \mathbb R^5$, 
\begin{equation*}
	\mathcal R_T(\tau_1,\tau_2)	=\bigg(\bigg(\int_{F_j}\tau\cdot n_{T}\ds \bigg)_{1\leq j\leq3 } ,\,
	[\tau\cdot n_\Gamma],\,
	%div\sigma_1 -div\sigma_2,\,
	[K^{-1}\tau\cdot t_\Gamma] (x_{\Gamma})\bigg).
\end{equation*}
\begin{lmm}\label{lem:bij}
The operator $\mathcal R_T$ is bijective. %due to the unisolvence of the $\IRT^0(T)$ space, and therefore surjective.	
\end{lmm}
\begin{proof}
Since $\mathrm{dim} \, \mathcal E(T)=5$, it is sufficient to prove that $\mathcal R_T$ is injective. 
Let $\tau=(\tau_1,\tau_2) \in \mathcal E(T)$ such that $\mathcal R_T(\tau_1,\tau_2)=0$, i.e.
\begin{equation*}
{\tau}_{|T^i} \in \RT^0(T^i),\quad 
\div(\tau)_{|T^1}=\div(\tau)_{|T^2},\quad 
[\tau\cdot n_\Gamma]=0,\quad 
[K^{-1}\tau\cdot t_\Gamma](x_\Gamma)=0
\end{equation*} 
and $N_{T,j}(\tau)=0$ for $1\leq j\leq 3$. But this implies that $\tau\in\IRT^0(T)$ and the unisolvency of this space immediately yields that $\tau=0$, \Red{whence} the result.
\end{proof}

\subsection{Construction of specific shape functions}\label{subsec:construction}

Since $\mathcal R_T$ is bijective, there exist unique functions $\Lambda_n$ and $\Lambda_t$ in $\mathcal E(T)$ such that 
\begin{equation}\label{eq: Lambda_n,t_vect}
	\mathcal R_T(\Lambda_n)=(0,0,0,1,0), \quad \mathcal R_T(\Lambda_t)=(0,0,0,0,1).	
\end{equation}
In other words, $\Lambda_n$ and $\Lambda_t$ belong to $ \mathcal{E}(T)$ and satisfy $N_{T,j}(\Lambda_t)=N_{T,j}(\Lambda_n)=0$ for $1\le j\le 3$, as well as: 
\begin{equation}\label{eq: ppt phi_n and phi_t}
	\left\{	\begin{array}{rl}
		& {\Lambda_n}_{|T^i} \in \RT^0(T^i),\quad i=1,2\\
		&\div({\Lambda_n}_{|T^1})=\div({\Lambda_n}_{|T^2})\\
		&[\Lambda_n\cdot n_\Gamma]=1\\
		&[K^{-1}\Lambda_n\cdot t_\Gamma](x_\Gamma)=0
	\end{array} \right. \,,
	\quad 
	\left\{\begin{array}{rl}
		&{\Lambda_t}_{|T^i} \in \RT^0(T^i),\quad i=1,2\\
		&\div({\Lambda_t}_{|T^1})=\div({\Lambda_t}_{|T^2})\\
		&[\Lambda_t\cdot n_\Gamma]=0\\
		&[K^{-1}\Lambda_t\cdot t_\Gamma](x_\Gamma)=1
	\end{array}\right. \Red{.}
\end{equation}    

In the sequel, we express $\Lambda_n$ and $\Lambda_t$ in terms of the functions $\phi_n$ and $\phi_t$ previously introduced in \eqref{eq:phi_t_n}. The only difference between $\phi_n,\,\phi_t$ and $\Lambda_n, \,\Lambda_t$  is the presence of the diffusion coefficient $K$ in the jump of the tangential trace. 

\begin{lmm}\label{thm:construction_lambda}
Let $k_\triangle:= k_{|T^\triangle}$ and $k_\square:= k_{|T^\square}$. One has that 
\begin{equation}\label{eq: def Lambda_n,t_IRT}
\Lambda_n= \alpha_n \phi_t+\phi_n, \quad \Lambda_t= \beta_t\phi_t
\end{equation}
where 
\begin{equation}\label{coeff_alpha_beta_A}
\alpha_n=\frac{(1- \frac{k_\triangle}{k_\square}) (I_T\rho\cdot t_\Gamma^\triangle)(x_\Gamma)}{A},\quad  \beta_t = \frac{k_\triangle}{A},\quad A=1- (1- \frac{k_\triangle}{k_\square}) (I_T\omega\cdot t_\Gamma^\Delta)(x_\Gamma). 
\end{equation}
\end{lmm}
\begin{proof}
Since we already have the existence and uniqueness of $\Lambda_n$ and $\Lambda_t$, we can start from \eqref{eq: def Lambda_n,t_IRT} and show that we can uniquely determine the constants $\alpha_n$ and $\beta_t$ such that $\Lambda_n$ and $\Lambda_t$ thus defined satisfy the desired properties \eqref{eq: ppt phi_n and phi_t}. 
Thanks to the properties of $\phi_n$ and $\phi_t$ stated in Lemma \ref{lem:ppty_phi}, we clearly have $\Lambda_n,\,\Lambda_t \in \mathcal E(T)$ and
\begin{equation*}
N_{T,j}(\Lambda_n)=N_{T,\Blue{j}}(\Lambda_t)=0\quad (1\le j\le 3),\qquad 
[\Lambda_n\cdot n_\Gamma] =1,\qquad 
[\Lambda_t\cdot n_\Gamma]=0.
\end{equation*}	

Next, we look for $\alpha_n$ and $\beta_t$ such that $[K^{-1}\Lambda_n\cdot t_\Gamma](x_\Gamma)=0$ and $[K^{-1}\Lambda_t\cdot t_\Gamma](x_\Gamma)=1$. We first impose
\begin{equation}\label{eq: help_init}
0=[K^{-1}\Lambda_n\cdot t_\Gamma](x_\Gamma)=\alpha_n [K^{-1}\phi_t\cdot t_\Gamma](x_\Gamma) + [K^{-1}\phi_n\cdot t_\Gamma](x_\Gamma).
\end{equation}
A simple computation yields: 
\begin{equation*}
\begin{split}
[K^{-1}\phi_t\cdot t_\Gamma](x_\Gamma)&=  (K^{-1}{\omega\cdot t_\Gamma})_{|T^1}(x_\Gamma) - (K^{-1}{\omega\cdot t_\Gamma})_{|T^2}(x_\Gamma)- [K^{-1}] (I_T\omega\cdot t_\Gamma)(x_\Gamma),\\
[K^{-1}\phi_n\cdot t_\Gamma](x_\Gamma)&= -[K^{-1}] (I_T\rho\cdot t_\Gamma)(x_\Gamma).
\end{split}
\end{equation*}
Next, we distinguish between two cases, depending on which subdomain contains $T^\triangle$. We define the jump with respect to $n_\Gamma^\triangle$ by $[w]_\triangle= n_\Gamma\cdot n_\Gamma^\triangle[w]$. If $T^\triangle=T^1$, then $t_\Gamma=t_\Gamma^\triangle$, $[K^{-1}]=[K^{-1}]_\triangle=k_\triangle^{-1}-k_\square^{-1}$ and the previous relations become:
\begin{equation}\label{eq: final_help}
[K^{-1}\phi_t\cdot t_\Gamma]= k_\triangle^{-1} - [K^{-1}]_\triangle I_T\omega\cdot t_\Gamma^\triangle,\qquad 
[K^{-1}\phi_n\cdot t_\Gamma]= -[K^{-1}]_\triangle I_T\rho\cdot t_\Gamma^\triangle.
\end{equation}
If $T^\triangle=T^2$, then $t_\Gamma=-t_\Gamma^\triangle$, $[K^{-1}]=-[K^{-1}]_\triangle$ and we obtain the same formulae \eqref{eq: final_help}. From  \eqref{eq: help_init} and \eqref{eq: final_help}, we obtain:
\begin{equation*}
\alpha_n \bigg( 1- (1- \frac{k_\triangle}{k_\square}) (I_T\omega\cdot t_\Gamma^\Delta)(x_\Gamma)\bigg)=(1- \frac{k_\triangle}{k_\square}) (I_T\rho\cdot t_\Gamma^\Delta)(x_\Gamma),
\end{equation*}
which yields the desired formula (\ref{coeff_alpha_beta_A}) for $\alpha_n$ if the denominator $A$ does not vanish.

We next show that $A$ indeed has a strictly positive lower bound, thanks to \eqref{eq: bound_omega_2}. If $ k_\triangle \geq k_\square$, then using $1- \frac{k_\triangle}{k_\square} \le 0$ and $ (I_T\omega\cdot t_\Gamma^\triangle)(x_\Gamma) \ge 0$ we get $A\geq 1$. If $ k_\triangle < k_\square$, then we use $1- \frac{k_\triangle}{k_\square} >0$ and $(I_T\omega\cdot t_\Gamma^\triangle)(x_\Gamma)\leq 1$ to get   
$A\geq1-(1-\dfrac{k_\triangle}{k_\square}) =\dfrac{k_\triangle}{k_\square}$.

Finally, let us focus on the coefficient $\beta_t$. Similarly to \eqref{eq: help_init}, we  now impose that   
\begin{equation*}
1=[K^{-1}\Lambda_t\cdot t_\Gamma](x_\Gamma)=\beta_t [K^{-1}\phi_t\cdot t_\Gamma](x_\Gamma),
\end{equation*}
which thanks to the first equality of \eqref{eq: final_help} translates into $\beta_t A=k_\triangle$
and thus yields the announced result.
\end{proof}
	
\subsection{Bounds of the specific shape functions}
It is useful to bound the $L^2$-norm of $\Lambda_n$ and $\Lambda_t$ in order to subsequently bound the a posteriori error estimator $\eta_T$. 

\begin{lmm}\label{trm: bound Lambda_n,t}
One has that:
\begin{equation*}
\nrmL[T]{\Lambda_n}\lesssim \frac{k_{max}}{k_{min}}  \sum_{j=2}^3|F_j^\triangle|, \quad \nrmL[T]{\Lambda_t}\lesssim k_{max}\sum_{j=2}^3|F_j^\triangle|.
\end{equation*}
%where we recall that $F_j^\triangle=F_j\cap \partial T^\triangle$.
\end{lmm}
\begin{proof} 
From Lemma \ref{thm:construction_lambda}, we have on the one hand, 
\begin{equation}\label{eq:estim_lambda}
\|\Lambda_n\|_{T} \leq |\alpha_n|\nrmL[T]{\phi_t}+ \nrmL[T]{\phi_n},\qquad \|\Lambda_t\|_{T} \leq |\beta_t| \nrmL[T]{\phi_t}.
\end{equation}
On the other hand, we have \Red{also showed in the proof of Lemma \ref{thm:construction_lambda}} that $A\ge \min\{1,\displaystyle \frac{k_\triangle}{k_\square}\}$, which allows us to bound $\alpha_n$ and $\beta_t$ in the sequel \Red{from their expressions given in \eqref{coeff_alpha_beta_A}.} If $ k_{max}=k_\triangle$, then using \eqref{eq: bound_rho} and $A\ge 1$, \Red{we get that}
\begin{equation}\label{eq: beta and alpha bound}
| \alpha_n |\lesssim \bigg|1-\frac{k_{max}}{k_{min}}\bigg| \lesssim   \frac{k_{max}}{k_{min}},\qquad 0<\beta_t\leq k_{max},
\end{equation}
whereas if $ k_{max}=k_\square$, we get, using  again \eqref{eq: bound_rho} and $A\ge \displaystyle\frac{k_\triangle}{k_\square}$, the same bound as above. In conclusion, $\alpha_n$ and $\beta_t$ satisfy \eqref {eq: beta and alpha bound}.

Thanks to Lemma \ref{lem:bound phi_t_n} and \eqref{eq:estim_lambda}, we can write that 
\begin{equation*}
\|\Lambda_n\|_{T}\lesssim  (|\alpha_n|+1) \sum_{j=2}^3|F_j^\triangle|,
\qquad \|\Lambda_t\|_{T} \lesssim |\beta_t| \sum_{j=2}^3|F_j^\triangle|
\end{equation*}
which combined with \eqref{eq: beta and alpha bound} ends the proof. 
\end{proof}

\Blue{
\begin{remark}\label{rk:2}
In the case the triangle $T\in \T^\Gamma$ is cut by one of its vertices into two sub-triangles, one may choose either of the two sub-triangles as $T^{\Delta}$. The only changes in Section \ref{sec:decomp_RT} occur in the sum $\sum_{j=2}^3|F_j^{\Delta}|$, which is now replaced by the portion of the cut edge contained in $\partial T^\Delta$, denoted $|F^\Delta|$. 
\end{remark}
}

\section{Decomposition of the flux correction in the cut cells}\label{sec:bound_tau}

In this section, we write $\sigIRT -K\nabla_h u_h$ on a cut cell $T$ in terms of the shape functions $\Lambda_n$ and $\Lambda_t$, in view of its further estimation in the next section. 
%Recall that $\sigIRT -K\nabla_h u_h$ is piecewise $\RT^0$ on any $T\in\T^{\Gamma}$, due to the discontinuity of the diffusion coefficients $ K $ and of the approximate solution $ u_h$.

First, we introduce $\varphi_0\in \RT^0(T)$ uniquely defined by the conditions:
\begin{equation}\label{eq: def_phi_0}
N_{T,j}(\varphi_0)=N_{T,j}(\sigIRT - K\nabla_h u_h)=\frac{1}{|F_j|}\int_{F_j}(\sigIRT - K\nabla_h u_h)\cdot n_{T}\ds,\quad  1\le j\le 3. 
\end{equation}
Furthermore, since $\sigIRT \in \IRT^0(T)$ and $\div \varphi_0\in P^0(T)$, one has
\[
\div (\sigIRT -k_1\nabla u_{h,1}-\varphi_0)_{|T^1}= \div (\sigIRT -k_2\nabla u_{h,2}-\varphi_0)_{|T^2}.
\]
Thus, $(\sigIRT -K\nabla_h u_h-\varphi_0)\in \mathcal E(T)$. Moreover, by definition of $ \varphi_0 $, one also has 
% and since  $\mathcal R_T$ is bijective,   then   there exists a unique  constants $C_n$ and $C_t$ such that 
\[
\int_{F_j} (\sigIRT-K\nabla_h u_h-\varphi_0)\cdot n_{T}\ds =0, \quad 1\leq j\leq 3.
\]
It follows that there exist $C_n, C_t \in \mathbb R$ such that $\mathcal R_T(\sigIRT -K\nabla_h u_h-\varphi_0)=(0,0,0,C_n,C_t)$. Since $\mathcal R_T$ is bijective according to Lemma \ref{lem:bij}, the previous relation combined with \eqref{eq: Lambda_n,t_vect} yields the following decomposition on $T$:
\begin{equation} \label{eq: decomposition}
\sigIRT - K\nabla_h u_h =\varphi_{0}+ C_n\Lambda_n+ C_t\Lambda_t.
\end{equation}

By writing that  
\begin{equation*}
\begin{split}
[ (\sigIRT-K\nabla_{\Blue{h}} u_h)\cdot n_\Gamma]=&[\varphi_{0}\cdot n_\Gamma]+C_n[\Lambda_n\cdot n_\Gamma]+ C_t[\Lambda_t\cdot n_\Gamma],\\
[K^{-1}(\sigIRT-\Blue{K}\nabla_{\Blue{h}} u_h)\cdot t_\Gamma] =& [K^{-1}]\varphi_{0}\cdot t_\Gamma+C_n[K^{-1}\Lambda_n\cdot t_\Gamma]+ C_t[K^{-1}\Lambda_t\cdot t_\Gamma]
\end{split}
\end{equation*}
and using the properties of $\varphi_0$, $\Lambda_n$, $\Lambda_t$ and $\sigIRT$, we can compute these constants:
\begin{equation}\label{eq: C_n,C_t}
C_n= - [K\nabla_{\Blue{h}} u_h\cdot n_\Gamma], \qquad  C_t= -[K^{-1}]\varphi_{0}\cdot t_\Gamma (x_\Gamma)- [\nabla\textbf{} u_h\cdot t_\Gamma](x_\Gamma).
\end{equation}

In the sequel, we are interested in bounding $\Vert \varphi_0\Vert_T$, $|C_n|$ and $|C_t|$. Clearly, \Blue{from \eqref{eq: bound base RT} we have}
\begin{equation}\label{eq: bound phi_00}
\nrmL[T]{\varphi_0}\leq \sum_{j=1}^3|N_{T,j}(\varphi_0)| \nrmL[T]{\Lambda_{T,j}}\lesssim h_T\sum_{j=1}^3|N_{T,j}(\varphi_0)|.
\end{equation}
Thanks to the definition of the flux $\sigIRT$, we get, for any edge $F_j\subset \partial T$, that
\begin{equation}\label{eq: def_phi_non_cut}
(n_T\cdot n_{F_j}) N_{T,j}(\varphi_0)=\frac{1}{|F_j|}\bigg(\Red{-\frac{n_T\cdot n_{F_j}}{2}}\int_{F_j}k_i\jF{\partial_n u_{h,i}}\ds \Red{-}\int_{F_j}k_i\theta_{h,i}\ds\bigg) \quad \text{if } F_j\in \F^i\setminus\F^{\Gamma},
\end{equation}
\begin{equation}\label{eq: def_phi_cut}
(n_T\cdot n_{F_j})N_{T,j}(\varphi_0)=\frac{1}{|F_j|}\sum_{i=1}^2\bigg(\Red{-\frac{n_T\cdot n_{F_j}}{2}}\int_{F_j^i}k_i\jF{\partial_n u_{h,i}}\ds \Red{-}\int_{F_j}k_i\theta_{h,i}\ds\bigg) \quad \text{if } F_j\in \F^{\Gamma}. 
\end{equation}

\begin{lmm}\label{lem: bound of N_i,T(phi_0)}
The following bound holds true:  
\begin{equation*}
\begin{split}
h_T\sum_{j=1}^3  |N_{T,j}(\varphi_0)|\lesssim &   \sum_{\Tp\in \Delta_T}\bigg(h_\Tp\|f\|_{\Tp}+h_\Tp^{1/2}\sum_{i=1}^2\|k_i\jF{\partial_nu_{h,i}}\|_{\partial \Tp\setminus \partial \Delta_T}\bigg)
+\sum_{i=1}^2 k_i^{1/2}j_{i,\Delta_T}(u_{h,i},u_{h,i})^{1/2}\\
 &+\sum_{\Tp\in\T^\Gamma\cap \Delta_T}\bigg(\frac{k_{min}}{h_\Tp^{1/2}}\|[u_h]\|_{\Gamma_\Tp}+h_\Tp^{1/2}\|[K\nabla u_h\cdot n_\Gamma]\|_{\Gamma_\Tp}\bigg).
\end{split}
\end{equation*}
\end{lmm}
\begin{proof}
Let $F_j$ a cut edge of $T$. Thanks to \eqref{eq: def_phi_cut}, the Cauchy-Schwarz inequality and $ |F_j^i|\leq |F_j|$, we have that: 
\begin{equation*}
|N_{T,j}(\varphi_0)|\leq\sum_{i=1}^2\bigg(\frac{1}{\Red{2}|F_j|^{1/2}}\nrmL[F_j]{k_i\jF{\partial_n u_{h,i}}}+ \frac{k_i}{|F_j|^{1/2}}\nrmL[F_j]{\theta_{h,i}}\bigg).
\end{equation*}
On a non-cut edge $F_j\in \F^i$, we obtain from  \eqref{eq: def_phi_non_cut} a similar bound:
\begin{equation*}
|N_{T,j}(\varphi_0)| \leq\frac{1}{\Red{2}|F_j|^{1/2}}\nrmL[F_j]{k_i \jF{\partial_n u_{h,i}}}+ \frac{k_i}{|F_j|^{1/2}}\nrmL[F_j]{\theta_{h,i}}.
\end{equation*}
The previous estimates lead to:
\begin{equation*}
h_T\sum_{j=1}^3  |N_{T,j}(\varphi_0)|\lesssim \sum_{i=1}^2 (h_T^{1/2}\nrmL[\partial T]{k_i\jF{\partial_n u_{h,i}}}+ k_i h_T^{1/2}  \nrmL[\partial T]{\theta_{h,i}}).
\end{equation*}
Next, the relation $(\theta_{h,i})_{|F}=(\theta_N^i)_{|F}+(\theta_M^i)_{|F}$ on any edge $F$ of vertices $N,M$, together with the bound of the local multiplier $\theta_N^i$ given in Lemma \ref{thm:bound_theta}, yields: 
\begin{equation*}
\begin{split}
\sum_{i=1}^2 k_ih_T^{1/2}\|\theta_{h,i}\|_{\partial T} \lesssim & \sum_{\Tp\in
\Delta_T} \left(h_\Tp\|f\|_{\Tp} +\sqrt{h_\Tp}  \sum_{i=1}^2 \| k_i \jF{\partial_n u_{h,i}}\|_{\partial \Tp\setminus \partial \Delta_T} \right)+\sum_{i=1}^2 k_i^{1/2}j_{i,\Delta_T}(u_{h,i},u_{h,i})^{1/2}\\
&+\sum_{\Tp\in\T^\Gamma\cap \Delta_T}
\left(\frac{{k_\Gamma}}{\sqrt{ h_\Tp}}\|\jump{u_h}\|_{\Gamma_\Tp}+
\sqrt{h_\Tp}\|[K\nabla u_h\cdot n_\Gamma]\|_{\Gamma_\Tp}\right).
\end{split}
\end{equation*}
Combining the two previous estimates, we get the announced result. 
\end{proof}
	
\begin{lmm} \label{lem: bound C_n,C_t}
One has that:
\begin{equation*}
|C_n|= \frac{1}{|\Gamma_T|^{1/2}}\nrmL[\Gamma_T]{[K\nabla u_h\cdot n_\Gamma]},\qquad
|C_t|\lesssim|[K^{-1}]|\sum_{j=1}^3|N_{T,j}(\varphi_{0})| + \frac{1}{|\Gamma_T|^{3/2}}\nrmL[\Gamma_T]{[u_h]}.
\end{equation*}
\end{lmm}

\begin{proof}
The first equality is obvious, it follows from \eqref{eq: C_n,C_t} using that $[K\nabla u_h\cdot n_\Gamma]$ is constant on $\Gamma_T$. From the definition \eqref{eq: C_n,C_t} of $C_t$, we get that:
\begin{equation}\label{eq: bound c_t_prf}
|C_t|\leq |[K^{-1}]||\varphi_{0}(x_\Gamma)\cdot t_\Gamma | + |[\nabla u_h\cdot t_\Gamma](x_\Gamma)|.
\end{equation}
			
On the one hand, we have thanks to \eqref{eq: bound base RT} that
\begin{equation} \label{eq: jump_2_phi_0}
|\varphi_{0}(x_\Gamma)\cdot t_\Gamma| \leq  \sum_{j=1}^3|N_{T,j}(\varphi_{0})||\Lambda_{T,j}(x_\Gamma)\cdot t_\Gamma|\lesssim \sum_{i=1}^3|N_{T,j}(\varphi_{0})|.
\end{equation}
			
On the other hand, by means of an inverse inequality\Red{,} we have:
%using that $[\nabla u_h\cdot t_\Gamma]$ is constant on $\Gamma_T$ and an inverse inequality for discrete functions, we get that:
\begin{equation}\label{eq: jump_final_phi_0}
|[\nabla u_h\cdot t_\Gamma](x_\Gamma)|^2 = \frac{1}{|\Gamma_T|}\int_{\Gamma_T}[\nabla u_h\cdot t_\Gamma]^2\ds\lesssim
\frac{1}{|\Gamma_T|^3} \int_{\Gamma_T}[u_h]^2\ds.
\end{equation}
			
We obtain the announced estimate for $C_t$ by using \eqref{eq: jump_2_phi_0} and \eqref{eq: jump_final_phi_0} in \eqref{eq: bound c_t_prf}.
\end{proof}

\section{Efficiency}\label{sec:eff} 
In this section, we prove the local efficiency of each indicator $\eta_T$, $\tilde \eta_T$ and $\eta_F$, with respect to the following norm of the error: 
\begin{equation*}%\label{eq: norm_h}
\|v_h\|_{h,\Delta_T}^2= \sum_{i=1}^{2}\bigg(k_i\|\nabla v_{h,i}\|^2_{\Delta_T\cap \Omega^i}+j_{i,\Delta_T}(v_{h,i},v_{h,i})\bigg)+\sum_{T\in \mathcal{T}_h^{\Gamma}\cap \Delta_T}\int_{\Gamma_T}\frac{k_{\Gamma}}{h_T}[v_h]^2 \ds \Red{.}
%+\sum_{T\in \mathcal{T}_h^{\Gamma}\cap \Delta_T}\int_{\Gamma_T}\frac{h_T}{k_{\Gamma}}\{K\nabla v_h\cdot n_\Gamma\}^2 \ds.
\end{equation*}
The main difficulty lies in the treatment of the cut elements, in particular as regards the dependence of the constants on the coefficients and the mesh/interface geometry.   

We begin with $\eta_T$ and we bound it by the local error, with a constant $C_{T}$ independent of the mesh, but dependent on the coefficients and the interface geometry and given explicitly. The proof for a cut triangle $ T $ uses the results of Sections \ref{sec:decomp_RT} and \ref{sec:bound_tau} and is quite technical; therefore, we present it in several steps.
	
\subsection{Bound of $\eta_T$ by the residuals}\label{subsec:eta_T1}
The first step consists \Red{of} bounding the local estimator $\eta_T$ by a residual-type estimator, for any $T\in \T$.

\begin{lmm}\label{lem: bound_1_effica}
If $T\in\T^i\backslash\T^\Gamma$ ($i=1,2$), one has that
\begin{equation}\label{eq:bound_noncut}
\eta_T\lesssim \sum_{T'\in\Delta_T}\bigg(\frac{h_{T'}}{k_{i}^{1/2}}\|f\|_{T'}+h_{T'}^{1/2}k_{i}^{1/2}(\|\jF{\partial_n u_{h,i}}\|_{\partial T'\setminus \partial \Delta_T}+\|\jF{\partial_n u_{h,i}}\|_{\partial T'\cap\Fg^i})\bigg)\Red{,}
\end{equation}
whereas if $T\in \T^\Gamma$, then  
\begin{equation}\label{eq: bound_cut_rsd}
\begin{split} 
\eta_T\lesssim& \frac{k_{max}}{k_{min}}\sum_{T'\in\Delta_T}\bigg(\frac{h_{T'}}{k_{min}^{1/2}}\|f\|_{T'}+\frac{h_{T'}^{1/2}}{k_{min}^{1/2}}\sum_{i=1}^2 \|k_i\jF{\partial_n u_{h,i}}\|_{\partial T'\setminus \partial \Delta_T}\bigg)
+\frac{k_{max}^{3/2}}{k_{min}^{3/2}}\sum_{i=1}^2j_{i,\Delta_T}(u_{h,i},u_{h,i})^{1/2}\\
&+\frac{k_{max}}{k_{min}}\sum_{T'\in\T^\Gamma\cap \Delta_T}\bigg(\frac{k_{min}^{1/2}}{|\Gamma_{T'}|^{1/2}}\|[u_h]\|_{\Gamma_{T'}}+\frac{h_{T'}^{1/2}}{k_{min}^{1/2}}\|[K\nabla u_h\cdot n]\|_{\Gamma_{T'}}\bigg).
\end{split}
\end{equation}
\end{lmm}
\begin{proof}
The proof for a non-cut triangle is similar to the one given in \Blue{Lemma 6.2 of} \cite{Aimene}, since $\sigIRT-k_i\nabla u_{h,i} \in \Red{\RT}^0(T)$. \Blue{In our setting, the Dirichlet condition on $\partial \Omega$ is imposed strongly, the diffusion coefficient is constant on each $\Omega^i$, and we use linear finite elements, which simplifies the estimate of $\eta_T$ given in \cite{Aimene}. In comparison to \cite{Aimene}, we have an additional term on $\mathcal F_g^i$ in the right-hand side of (\ref{eq:bound_noncut}). This term arises from the bound of $\theta_{h,i}$ given in Lemma \ref{thm:bound_theta} and accounts for the ghost penalization, which was not present in \cite{Aimene}.}

So in the sequel, we only treat the case $T\in \T^\Gamma$. We have, thanks to \eqref{eq: decomposition} and to the triangle inequality, that
\begin{equation*}
\|\sigIRT-K\nabla_h u_h\|_{T}\leq  \|\varphi_{0}\|_{T}+|C_n|\| \Lambda_n\|_{T}+ |C_t|\| \Lambda_t\|_{T}.
\end{equation*}
From \eqref{eq: bound phi_00}, Lemma \ref{trm: bound Lambda_n,t} and Lemma \ref{lem: bound C_n,C_t}, we  get the following bound: 
\begin{equation*}
\begin{split}
\nrmL[T]{\sigIRT-K\nabla_h u_h} &\lesssim  h_T\sum_{j=1}^3|N_{T,j}(\varphi_0)|
+\frac{k_{max}}{k_{min}} (\sum_{j=2}^3|F_j^\triangle|)|\Gamma_T|^{-1/2}\nrmL[\Gamma_T]{[K\nabla u_h\cdot n_\Gamma]}\\
&+ k_{max}(\sum_{j=2}^3|F_j^\triangle|)\bigg( |[K^{-1}]|\sum_{j=1}^3|N_{T,j}(\varphi_{0})| + |\Gamma_T|^{-3/2}\nrmL[\Gamma_T]{[u_h]}\bigg). 
\end{split}
\end{equation*}
Using that $k_{max}|[K^{-1}]|=\frac{k_{max}}{k_{min}}-1$ and $|F_j^\triangle|\leq h_T$,  we next get: 
\begin{equation*}
\begin{split}
\nrmL[T]{\sigIRT-K\nabla_h u_h}\lesssim & \frac{k_{max}}{k_{min}}\bigg( h_T\sum_{j=1}^3|N_{T,j}(\varphi_0)|
+(\sum_{j=2}^3|F_j^\triangle|)|\Gamma_T|^{-1/2}\nrmL[\Gamma_T]{[K\nabla u_h\cdot n_\Gamma]}\\
&+ k_{min}(\sum_{j=2}^3|F_j^\triangle|)|\Gamma_T|^{-3/2}\nrmL[\Gamma_T]{[u_h]}  \bigg).
\end{split}
\end{equation*}

Note that $\displaystyle \sum_{j=2}^3|F_j^\triangle|\lesssim |\Gamma_T|$, even in the worst cut case where $|\Gamma_T|$ is small, since then $\displaystyle \sum_{j=2}^3|F_j^\triangle|$ is of the same size as $|\Gamma_T|$. Using in addition $|F_j^\triangle|\leq h_T$ and $ |\Gamma_T|\leq h_T$, we further get: 
\begin{equation*}
\eta_T\leq \frac{1}{k_{min}^{1/2}}  \nrmL[T]{\sigIRT-K\nabla_h u_h}
\lesssim \frac{k_{max}}{k_{min}}\bigg( \frac{h_T}{k_{min}^{1/2}} \sum_{j=1}^3|N_{T,j}(\varphi_0)|+  \frac{h_T^{1/2}}{k_{min}^{1/2}}\nrmL[\Gamma_T]{[K\nabla u_h\cdot n_\Gamma]}
+ \frac{k_{min}^{1/2}}{|\Gamma_T|^{1/2}}\nrmL[\Gamma_T]{[u_h]}\bigg).
\end{equation*}
We next apply Lemma \ref{lem: bound of N_i,T(phi_0)}; by using $h_T^{-1/2}\leq |\Gamma_T|^{-1/2}$ in front of $[u_h]$, we obtain \eqref{eq: bound_cut_rsd} which ends the proof.
\end{proof}

\Blue{
\begin{remark}
The proof of Lemma \ref{lem: bound_1_effica} remains valid if $T\in\T^{\Gamma}$ is cut by one of its vertices into two sub-triangles. In this case, $\sum_{j=2}^3|F_j^{\Delta}|$ is replaced by $|F^\Delta|$, where $F^\Delta$ is the portion of the cut edge contained in $\partial T^\Delta$ (cf. Remark \ref{rk:2}) and we still have $|F^\Delta|\lesssim |\Gamma_T|$.
\end{remark}
}

\subsection{Bound of the residuals by the local error}\label{subsec:eta_T_2}
The second step consists \Red{of} bounding the residual terms on the right-hand side of \eqref{eq:bound_noncut} and \eqref{eq: bound_cut_rsd} by the error for any $T\in \T$. First, let us note that thanks to the triangle inequality
\begin{equation*}
\nrmL[T]{f} \le \nrmL[T]{f_h}+\nrmL[T]{f-f_h},
\end{equation*}
it is sufficient to bound $\nrmL[T]{f_h}$ instead of $\nrmL[T]{f}$, the last term being contained in $\epsilon(T)$. \Blue{The definitions of $f_h$ and $\epsilon(T)$ can be found in Section \ref{sec: A posteriori_estim}.}

We begin by treating the jump of the normal trace across an edge $F$ of $T'\in \Delta_T$. 

\Blue{\begin{lmm}\label{lem:new_lemma}
For $T\in\T^i \setminus \T^{\Gamma}$ ($1\le i\le 2$), one has
\begin{equation}\label{eq:estim_resid_term1}
\eta_T\lesssim  |u-u_{h,i}|_{1,K,\Delta_T}+j_{i,\Delta_T}(u_{h,i},u_{h,i})^{1/2}+\epsilon(\Delta_T)+\sum_{T'\in\Delta_T}\frac{h_{T'}}{k_{i}^{1/2}}\|f_h\|_{T'},
\end{equation}
whereas for $T\in\T^{\Gamma}$, one has
\begin{equation}\label{eq:estim_resid_term2}
\begin{split} 
\eta_T\lesssim & \frac{k_{max}^{3/2}}{k_{min}^{3/2}}(|u-u_h|_{1,K,\Delta_T}+\sum_{i=1}^2j_{i,\Delta_T}(u_{h,i},u_{h,i})^{1/2})+\frac{k_{max}}{k_{min}^{3/2}}\sum_{T'\in\Delta_T} h_{T'}\|f_h\|_{T'}\\
&+\frac{k_{max}}{k_{min}}\sum_{T'\in\T^\Gamma\cap \Delta_T}\bigg(\frac{k_{min}^{1/2}}{|\Gamma_{T'}|^{1/2}}\|[u_h]\|_{\Gamma_{T'}}+\frac{h_{T'}^{1/2}}{k_{min}^{1/2}}\|[K\nabla u_h\cdot n]\|_{\Gamma_{T'}}\bigg)+\frac{k_{max}^{3/2}}{k_{min}^{3/2}}\epsilon(\Delta_T).
\end{split}
\end{equation}
\end{lmm}
}
\begin{proof}
Assume first $T\in\T^i \setminus \T^{\Gamma}$. \Blue{Let $T'\in\Delta_T$ and $F\subset \partial T'$. Thanks to the shape regularity of the triangle $T'$, we have that $h_{T'}\lesssim h_F$.} If $F\in \Fg^i$, then the term $h_{T'}^{1/2}k_i^{1/2}\|\jF{\partial_n u_{h,i}}\|_{F}$ is easily bounded by the ghost penalty term $j_i(u_{h,i},u_{h,i})^{1/2}$ restricted to $F$. Otherwise, we have $ T^+_F\cup T_F^-\subset \Omega^i$ and the well-known Verf\"urth argument (cf. \cite{verfurth1994posteriori}) yields:
\begin{equation*}
    h_{T'}^{1/2}k_i^{1/2}\|\jF{\partial_n u_{h,i}}\|_{F}\lesssim |u-u_{h,i}|_{1,K,T^+_F\cup T_F^-}+ \epsilon(T^+_F\cup T_F^-).
\end{equation*}
For $F\subset \partial T'\setminus \partial \Delta_T $ we also have $T^+_F\cup T_F^-\subset \Delta_T$, so
(\ref{eq:estim_resid_term1}) follows from \eqref{eq:bound_noncut}.

Assume now $T\in\T^{\Gamma}$ \Blue{and let again $T'\in\Delta_T$}. For $F\subset \partial T'\setminus \partial \Delta_T $, we bound for $i=1,2$ the term $h_{T'}^{1/2}k_i^{1/2}\|\jF{\partial_n u_{h,i}}\|_{F}$ as previously, by means of Verf\"urth' argument. \Blue{Using furthermore that $\|k_i\jF{\partial_n u_{h,i}}\|_{F}\le k_{max}^{1/2}k_i^{1/2}\|\jF{\partial_n u_{h,i}}\|_{F}$,} we finally obtain (\ref{eq:estim_resid_term2}) from \eqref{eq: bound_cut_rsd}.
\end{proof}

Next, we focus on the term containing the data $f_h$.

\begin{lmm}\label{lem: bound_2_effica}
If $T\in \T^i\setminus\T^\Gamma$, then
\begin{equation*}
h_T\nrmL[T]{f_h} \lesssim  \nrmL[T]{k_i\nabla (u-u_{h,i})} +h_T\nrmL[T]{f-f_h},
\end{equation*}
whereas if $T\in \T^\Gamma$, one has:
\begin{equation*}
h_T\nrmL[T]{f_h} \lesssim \sum_{i=1}^{2} \nrmL[T^i]{k_i\nabla (u-u_{h,i})}+ h_T^{1/2} \nrmL[\Gamma_T]{ [ K \nabla u_h\cdot n_\Gamma]} +h_T\nrmL[T]{f-f_h}.
\end{equation*}
\end{lmm}
\begin{proof}
We only treat below the case of a cut triangle, the other one being standard.  Following Verf\"urth's argument (\Blue{see, for instance, Proposition 4.2. of} \cite{verfurth1994posteriori}), we let $b_T=27\lambda_1\lambda_2\lambda_3$ be the cubic bubble function on $T$ and consider  $w_T=f_hb_T$, which satisfies:
\begin{equation}\label{eq: ppt bubl_func}
\nrmL[T]{f_h}^2\lesssim \int_T f_hw_T\dx, \qquad
\nrmL[T]{\nabla w_T}\lesssim h_T^{-1} \nrmL[T]{ w_T}\lesssim  h_T^{-1} \nrmL[T]{f_h}.
\end{equation}
We start with 
\begin{equation*}
\nrmL[T]{f_h}^2\lesssim \int_T fw_T\dx +\int_T(f_h-f)w_T\dx
\lesssim  \int_T fw_T\dx +\nrmL[T]{f-f_h} \nrmL[T]{f_h}.
\end{equation*}
\Blue{Since $k_i\nabla u_{h,i}$ is piecewise constant on $\Omega_h^i$ for $i=1,2$, it follows that} $\div (k_i\nabla u_{h,i})_{|T^i}=0$. Using also that ${w_T}_{|\partial T}=0$, we get by integrating by parts on each $T^i=T\cap \Omega^i$: 
\begin{equation*}
\begin{split}
\int_T fw_T\dx &=- \sum_{i=1}^{2}  \int_{T^i} \div(k_i \nabla (u-u_{h,i}))w_T\dx \\
&=\sum_{i=1}^{2} \int_{T^i} k_i \nabla (u-u_{h,i})\cdot \nabla w_T\dx -\int_{\Gamma_T} [ K \nabla (u-u_h)\cdot n_\Gamma] w_T\ds\\
&\lesssim \nrmL[T]{\nabla w_T}(\sum_{i=1}^{2} \nrmL[T^i]{k_i\nabla (u-u_{h,i})})+ \nrmL[\Gamma_T]{w_T} \nrmL[\Gamma_T]{ [ K \nabla (u-u_h)\cdot n_\Gamma]}.
\end{split}
\end{equation*}
From the trace inequality \eqref{eq: trace_in_cut} and the properties \eqref{eq: ppt bubl_func} of $w_T$, it follows that:
\begin{equation*}
\nrmL[\Gamma_T]{w_T}\lesssim \frac{1}{h_T^{1/2}}  \nrmL[T]{w_T}+ h_T^{1/2}\nrmL[T]{\nabla w_T}\lesssim  h_T^{-1/2} \nrmL[T]{f_h},
\end{equation*}
which finally implies, thanks to $[ K \nabla u\cdot n_\Gamma]=0 $, that
\begin{equation*}
\nrmL[T]{f_h} \lesssim h_T^{-1}\sum_{i=1}^{2} \nrmL[T^i]{k_i\nabla (u-u_{h,i})}+ h_T^{-1/2} \nrmL[\Gamma_T]{ [ K \nabla u_h\cdot n_\Gamma]}+\nrmL[T]{f-f_h}.
\end{equation*}
Multiplying the previous inequality by $h_T$ yields the announced result.
\end{proof}
	
Finally, we are interested in bounding the term $\|[ K \nabla u_h\cdot n_\Gamma]\|_{\Gamma_T}$. For theoretical reasons only, to obtain the best constant, we make the following assumption.

\begin{assumption}\label{ass: 1}
For any $T\in \T^\Gamma$, there exist closed, regular shaped triangles $\tilde{T}^1\subset \Omega^1,\tilde{T}^2\subset \Omega^2$ 
such that  they have  $\Gamma_T$ as a common side: $ \tilde{T}^1\cap \tilde{T}^2=\Gamma_T$.
\end{assumption}

Since $|\Gamma_T|\le h_T$, we can moreover suppose, without loss of generality, that $\tilde{T}^i$ are included in the neighborhood of $T$ consisting of the triangles sharing an edge with $T$.
\begin{remark}
The particular situation where Assumption 1 might not hold is discussed \Magenta{in \cite{Aimene_these}, page 156}.
%Appendix \ref{sec:appendixB}. 
It leads to a similar estimate, with an efficiency constant multiplied by $k_{max}^{1/2}/k_{min}^{1/2}$. However, we have not noticed any influence of Assumption 1 in the numerical experiments, see \cite{Article1}.
\end{remark}

\begin{lmm}\label{lem: bound_3_effica}
Under Assumption \ref{ass: 1}, we have for any $T\in \T^\Gamma$ that: 
\begin{equation*}%\label{eq:lem_bound_3_effica}
h_T^{1/2}  \|[K\nabla u_h\cdot n_\Gamma]\|_{\Gamma_T}\lesssim \sum_{i=1}^2 \bigg(\frac{h_T^{1/2} k_i^{1/2}}{|\Gamma_T|^{1/2}} \nrmL[\tilde T^i]{k_i^{1/2}\nabla (u_{h,i}- u)} + h_T\nrmL[\tilde T^i]{f-f_h}\bigg)
 + \sum_{i=1}^2 k_i^{1/2} j_{i,\Delta_T}(u_{h,i},u_{h,i})^{1/2}.
\end{equation*}
\end{lmm}
\begin{proof}
We use a similar argument to \cite{he2019residual} for the immersed finite element method, i.e. we consider the bubble function $b_\Gamma$ on $\Gamma_T$ which takes the value $1$ at the midpoint of $\Gamma_T$ and satisfies:
$(b_\Gamma)_{|\partial (\tilde{T^1}\cup \tilde{T^2})}=0$ and $(b_\Gamma)_{|\tilde{T}^i}\in P^2(\tilde{T^i})$ for $i=1,2$. \Blue{It is useful to note that for any regular-shaped triangle $T^*\in \{\tilde T^1,\tilde T^2\}$, we have:
$$ 
\|\nabla b_\Gamma\|_{T^*}\lesssim \frac{1}{|\Gamma_T|} \|b_\Gamma\|_{T^*}\lesssim 1.
$$
}
Let then $\psi_\Gamma=[K\nabla u_h\cdot n_\Gamma]b_\Gamma$, which satisfies, for $T^*\in \{\tilde T^1,\tilde T^2\}$,
\begin{equation}\label{eq: psi_bubl_proprerties}
\nrmL[T^*]{\nabla \psi_\Gamma}\lesssim \frac{1}{|\Gamma_T|} \nrmL[T^*]{\psi_\Gamma}\lesssim \frac{1}{|\Gamma_T|^{1/2}}  \|[K\nabla u_h\cdot n_\Gamma]\|_{\Gamma_T},
\end{equation}
\begin{equation}\label{eq: psi_bubl_proprerties2}
\|[K\nabla u_h\cdot n_\Gamma]\|_{\Gamma_T}^2\lesssim \bigg|\int_{\Gamma_T}[K\nabla u_h\cdot n_\Gamma]\psi_\Gamma\ds \bigg|.
\end{equation}
Using the continuity of $K\nabla u\cdot n_\Gamma$ across $\Gamma$ and integrating by parts on each $\tilde T^i$, we get, since  $(\psi_\Gamma)_{|(\partial \tilde T^i\setminus \Gamma_T)}=0 $,
\begin{equation*}
\int_{\Gamma_T}[K\nabla u_h\cdot n_\Gamma]\psi_\Gamma\ds =  \sum_{i=1}^2 \bigg(\int_{\tilde{T}^i}   k_i\nabla (u_{h,i}- u)\cdot \nabla \psi_\Gamma\dx +\int_{\tilde{T}^i} f\psi_\Gamma \dx
-\sum_{F^*\in \Blue{\mathcal F(\tilde T^i)}}\int_{F^*} k_i\jF{\partial_n(u_{h,i}-u)}\psi_\Gamma\ds\bigg),
\end{equation*}
\Blue{where the set $\mathcal F(\tilde T^i)$ contains the possible intersections of the sides of $\mathcal F_h^i$ with $\tilde T^i$.} Using next the Cauchy-Schwarz inequality, we get:
\begin{equation}\label{eq: help__1}
\begin{split}
\bigg| \int_{\Gamma_T}[K\nabla u_h\cdot n_\Gamma]\psi_\Gamma\ds \bigg|
\le  &\sum_{i=1}^2 (\nrmL[\tilde{T}^i]{k_i\nabla (u_{h,i}- u)}\nrmL[\tilde{T}^i]{\nabla \psi_\Gamma} + \nrmL[\tilde{T}^i]{f}\nrmL[\tilde{T}^i]{\psi_\Gamma})\\
&+\sum_{i=1}^2\sum_{F^*\in \Blue{\mathcal F(\tilde T^i)}} \nrmL[F^*]{k_i\jF{\partial_n u_{h,i}}} \nrmL[F^*]{\psi_\Gamma}.
\end{split}
\end{equation}

Let $F^*\in \Blue{\mathcal F(\tilde T^i)}$. We next apply the trace inequality \eqref{eq: trace_in_cut} to the discrete function $\psi_{\Gamma}$ on the triangle $\tilde T^i$ cut by $F^*$; using also an inverse inequality and the fact that $\tilde T^i$ is a regular triangle of size $|\Gamma_T|$, we obtain: 
$\nrmL[F^*]{\psi_\Gamma}\lesssim |\Gamma_T|^{-1/2}\nrmL[\tilde T^i]{\psi_\Gamma}$. Together with the properties \eqref{eq: psi_bubl_proprerties} and \eqref{eq: psi_bubl_proprerties2} of $\psi_\Gamma$, the inequality \eqref{eq: help__1} yields:
\begin{equation}\label{eq:estim_help}
\|[K\nabla u_h\cdot n_\Gamma]\|_{\Gamma_T}\lesssim \sum_{i=1}^2 \bigg(\frac{k_i^{1/2}}{|\Gamma_T|^{1/2}} \nrmL[\tilde T^i]{k_i^{1/2}\nabla (u_{h,i}- u)}+ |\Gamma_T|^{1/2} \nrmL[\tilde T^i]{f}+\sum_{F^*\in \Blue{\mathcal F(\tilde T^i)}} \nrmL[F^*]{k_i\jF{\partial_n u_{h,i}}}\bigg).
\end{equation}

Thanks to Verf\"urth's argument with the cubic bubble function on the triangle $\tilde T^i$ (of size $|\Gamma_T|$) we have, as in Lemma \ref{lem: bound_2_effica}, the following estimate: 
\begin{equation*}
\nrmL[\tilde T^i]{f_h} \lesssim \frac{1}{|\Gamma_T|}\nrmL[\tilde T^i]{k_{i}\nabla ( u_{h,i}-u)} \Red{+\frac{1}{|\Gamma_T|^{1/2}}\sum_{F^*\in \mathcal F(\tilde T^i)} \nrmL[F^*]{k_i\jF{\partial_n u_{h,i}}}}+\nrmL[\tilde T^i]{f-f_h},\qquad i=1, 2.
\end{equation*}
Using also the triangle inequality and $|\Gamma_T|\le h_T$, we finally get, for $i=1,2$,
\begin{equation}\label{eq:estim_f_help}
|\Gamma_T|^{1/2}\nrmL[\tilde T^i]{f} \lesssim \frac{k_{i}^{1/2}}{|\Gamma_T|^{1/2}}\nrmL[\tilde T^i]{k_{i}^{1/2}\nabla ( u_{h,i} -u)} 
\Red{+\sum_{F^*\in \mathcal F(\tilde T^i)} \nrmL[F^*]{k_i\jF{\partial_n u_{h,i}}}}
+h_T^{1/2}\nrmL[\tilde T^i]{f-f_h}.
\end{equation}

In addition, we note that for $i=1, 2$, \Blue{any $F^*\in \mathcal F_h(\tilde T^i)$ is included in a side of $\mathcal F_g^i\cap \mathcal F_{T'}$ for a triangle $T'\in\Delta_T$}. \Blue{Together with $h_T\lesssim h_F$}, this gives:
\begin{equation}\label{eq:estim_help_3}
\sum_{F^*\in \Blue{\mathcal F(\tilde T^i)}} \nrmL[F^*]{k_i\jF{\partial_n u_{h,i}}}\lesssim \Blue{\sum_{T'\in \Delta_T}}\sum_{F\in \mathcal F_g^i\cap \Blue{\mathcal F_{T'}}} \nrmL[F]{k_i\jF{\partial_n u_{h,i}}}\lesssim \frac{k_i^{1/2}}{h_T^{1/2}}j_{i,\Delta_T}(u_{h,i},u_{h,i})^{1/2}.
\end{equation}

We conclude using \eqref{eq:estim_f_help} and \eqref{eq:estim_help_3} in \eqref{eq:estim_help} and multiplying by $h_T^{1/2}$.
\end{proof}

\subsection{Final bound of $\eta_T$}
We now give the main result of this section, which ascertains the local efficiency of the error estimator $\eta_T$.

\begin{thrm}\label{thm: efficiency_IRT}
Under Assumption \ref{ass: 1}, for any $ T\in \T$ there exists a positive constant $C_T$ such that
\begin{equation*}
		\forall T\in \T,\quad \eta_T \lesssim C_{T} \left(\|u-u_h\|_{h,\Delta_T}+\epsilon(\Delta_T)\right),
\end{equation*}
with $C_T=1$ if $T\in \T\backslash \T^\Gamma$ and $C_{T}= \displaystyle \max_{T'\in \Delta_T \cap \T^{\Gamma}}\frac{h_{T'}^{1/2}}{|\Gamma_{T'}|^{1/2}}\frac{k_{max}^{3/2}}{k_{min}^{3/2}}$ if $T\in \T^\Gamma$.
\end{thrm}

\begin{proof}
It is sufficient to gather the results that we have proved so far in order to bound $\eta_T$. The estimate for a non-cut triangle is obvious, it follows from \eqref{eq:estim_resid_term1} combined with the first estimate of Lemma \ref{lem: bound_2_effica}. So we assume next that $T\in \T^\Gamma$. 

%Noting that $\tilde T^1\cup \tilde T^2 \subset \Delta_T$, 
Thanks to Lemma \ref{lem: bound_3_effica}, we can improve the second estimate of Lemma \ref{lem: bound_2_effica}: 
\begin{equation*}
\begin{split}
\frac{h_T}{k_{min}^{1/2}}\nrmL[T]{f_h}\lesssim 
& \frac{h_T}{k_{min}^{1/2}}(\nrmL[T]{f-f_h}+\nrmL[\tilde T^1\cup \tilde T^2 ]{f-f_h})\\
&+\frac{k_{max}^{1/2}}{k_{min}^{1/2}}\sum_{i=1}^{2}(|u-u_{h,i}|_{1,K,T^i} +\frac{h_T^{1/2}}{|\Gamma_T|^{1/2}} |u-u_{h,i}|_{1,K,\tilde T^i}+j_{i,\Delta_T}(u_{h,i},u_{h,i})^{1/2}).
\end{split}
\end{equation*}
We use this estimate on any cut triangle $T'\in\T^\Gamma\cap \Delta_T$, as well as Lemma \ref{lem: bound_3_effica} in the inequality \eqref{eq:estim_resid_term2}. \Blue{With $h_{T^{\prime}} \geq\left|\Gamma_{T^{\prime}}\right|$,} we finally obtain:
\begin{equation*}
\begin{split}
\eta_T\lesssim &\frac{k_{max}^{3/2}}{k_{min}^{3/2}}\bigg(\max_{T'\in \Delta_T \cap \T^{\Gamma}}\frac{h_{T'}^{1/2}}{|\Gamma_{T'}|^{1/2}} |u-u_h|_{1,K,\Delta_T}+\sum_{i=1}^2j_{i,\Delta_T}(u_{h,i},u_{h,i})^{1/2}+\epsilon(\Delta_T)\bigg)\\
&+\frac{k_{max}}{k_{min}}\sum_{T'\in\T^\Gamma\cap \Delta_T}\frac{h_{T'}^{1/2}}{|\Gamma_{T'}|^{1/2}} \frac{k_{\Gamma}^{1/2}}{h_{T'}^{1/2}}\|[u_h]\|_{\Gamma_{T'}}
\lesssim  C_T \|u-u_h\|_{h,\Delta_T}+\frac{k_{max}^{3/2}}{k_{min}^{3/2}}\epsilon(\Delta_T)
\end{split}
\end{equation*}
which concludes the proof.
\end{proof}

\subsection{Bounds of $\tilde \eta_T$ and $ \eta_F$}

The bound of $\tilde\eta_T$ on any cut triangle $ T\in \T^\Gamma $ is immediate, since $\tilde\eta_T$ already appears in the error $\|\cdot\|_{h,T}$, with a different scaling.
\begin{thrm}\label{lem: eta_tilde}
Let $ T\in \T^\Gamma $. There exists a positive constant $\tilde C_T$  such that
\begin{equation*}
\tilde \eta_T \le \tilde C_T\|u-u_h\|_{h,T},\quad \text{with}\quad \tilde C_T= \frac{h_T}{\sqrt{h_T^{min}|\Gamma_T|}}.
\end{equation*}
\end{thrm}
\begin{proof}
By definition, we have that: 
\begin{equation*}
\tilde \eta_T =\frac{\sqrt{h_T k_\Gamma}}{\sqrt{h_T^{min}|\Gamma_T|}}\nrmL[\Gamma_T]{[u_h]} = \frac{h_T}{\sqrt{\Red{h_T^{min}|\Gamma_T|}}}\frac{k_\Gamma^{1/2}}{h_T^{1/2}}\nrmL[\Gamma_T]{[u_h]}\le \tilde C_T\|u-u_h\|_{h,T},
\end{equation*}
which is the desired estimate.
\end{proof}
     		
It remains to bound the \Red{indicator} $\eta_F$ on a cut edge.
\begin{thrm}\label{thrm: bound_1_F}
Let $F\in\F^\Gamma$. There exists a positive constant $C_F$ such that
\begin{equation*}
\eta_F \lesssim C_{F}  \left(\|u-u_h\|_{h,\Delta_F}+\epsilon(\Delta_F)\right)\,\,\text{with }\,\, C_F=\max_{T\in \Delta_F\cap \T^{\Gamma}}\frac{h_T}{|\Gamma_T|}\frac{k_{max}^{3/2}}{k_{min}^{3/2}}.
\end{equation*}
\end{thrm}
\begin{proof}
By the definition of $\eta_F$ and the inequality \Blue{$\|a+b\|^2\le 2(\|a\|^2+\|b\|^2)$}, we have that:
\begin{equation*}
\begin{split}
\eta_F^2
&\Blue{=\frac{h_F}{k_{\Gamma}}\sum_{i=1}^2\int_{F^i} \jF{(\sigIRT-k_i\nabla u_{h,i})\cdot n_F+k_i\nabla u_{h,i}\cdot n_F }^2\ds}\\
&\le \frac{2h_F}{k_{\Gamma}}\sum_{i=1}^2\bigg( \int_{F^i} \jF{\sigma_h\cdot n_F-k_i\nabla u_{h,i}\cdot n_F }^2\ds+ \int_{F^i}k_i^2\jF{\partial_n u_{h,i}}^2\ds \bigg).
\end{split}
\end{equation*}
Using for the last integral that $F^i\subset F$, we get that 
\begin{equation}\label{eq:eta_F_1}
\eta_F^2 \le \frac{2h_F}{k_{\Gamma}}\nrmL[F] {\jF{(\sigIRT - K\nabla_h u_h)\cdot n_F }}^2+  \sum_{i=1}^2\frac{2k_i}{k_\Gamma}j_{i,\Delta_F}(u_{h,i},u_{h,i}).
\end{equation}
Since $k_i/k_\Gamma\lesssim k_{max}/k_{min}\le C_F\le C_F^2$, we only have to bound the first term on the right-hand side of (\ref{eq:eta_F_1}) to establish the theorem. 

Let $\{F\}=\partial T_F^+\cap \partial T_F^-$. By the triangle inequality, it is sufficient to bound $\nrmL[F] {(\sigIRT - K\nabla_h u_h)_{|T}\cdot n_F }$ for any $T\in \{T_F^+,T_F^-\}$. Using the decomposition \eqref{eq: decomposition} of  $\sigIRT - K\nabla_h u_h$ on $T$, we get: 
\begin{equation}\label{eq:eta_F_2}
\nrmL[F] {(\sigIRT - K\nabla_h u_h)_{|T}\cdot n_F }\leq \nrmL[F]{\varphi_0\cdot n_F}+ |C_n| 	\nrmL[F]{\Lambda_n\cdot n_F}+|C_t| 	\nrmL[F]{\Lambda_t\cdot n_F}.
\end{equation}
%where $\varphi_0$, $\Lambda_n$ and $\Lambda_t$ are given in  \eqref{eq: def_phi_0} and \eqref{eq: Lambda_n,t_vect}, respectively.			
For  $\varphi_0\in \RT^0(T)$, it is known that
\begin{equation}\label{eq: bound_phi_0_F}
\sqrt{h_F}\nrmL[F]{\varphi_0\cdot n_F}\lesssim \nrmL[T]{\varphi_0}\leq h_T\sum_{j=1}^3|N_{T,j}(\varphi_0)|,
\end{equation}
whereas from the definition \eqref{eq: def Lambda_n,t_IRT} of $\Lambda_n$ and $\Lambda_t$, together with the property $ (I_T w)\cdot n_F= \pi_F^0(w\cdot n_F)$, we get: 
\begin{equation*}
\nrmL[F]{\Lambda_n\cdot n_F}\le 2|\alpha_n| \nrmL[F]{\omega\cdot n_F}+ 2\nrmL[F]{\rho\cdot n_F},\qquad \nrmL[F]{\Lambda_t\cdot n_F}\le 2|\beta_t|\nrmL[F]{\rho\cdot n_F}.
\end{equation*}  			
Thanks to the definition (\ref{eq:omega_rho}) of $\omega$ and $\rho$, it follows, \Red{with $F^{\Delta}=F\cap \partial T^{\Delta}$,} that:
\begin{equation*}
\nrmL[F]{\rho\cdot n_F}^2= \int_{F^\triangle} |\
n_\Gamma^\triangle\cdot n_F|^2\ds \leq |F^\triangle|,\quad
\nrmL[F]{\omega\cdot n_F}^2= \int_{F^\triangle} |\
t_\Gamma^\triangle\cdot n_F|^2\ds \leq |F^\triangle|.
\end{equation*}
Using the bound (\ref {eq: beta and alpha bound}) for the coefficients $\alpha_n$ and $\beta_t$ and Lemma \ref{lem: bound C_n,C_t}, we end up with:
%\begin{equation*}
%\nrmL[F]{\Lambda_n\cdot n_F}\lesssim \sqrt{ |F^\triangle|} \frac{k_{max}}{k_{min}},\qquad 
%\nrmL[F]{\Lambda_t\cdot n_F}\lesssim \sqrt{|F^\triangle|} k_{max}.
%\end{equation*}    

\begin{eqnarray*}
|C_n|\nrmL[F]{\Lambda_n\cdot n_F}&\lesssim &\dfrac{k_{max}}{k_{min}} \dfrac{|F^\triangle|^{1/2}}{|\Gamma_T|^{1/2}}\nrmL[\Gamma_T]{[K\nabla u_h\cdot n_{\Gamma}]},\\
 |C_t|\nrmL[F]{\Lambda_t\cdot n_F}&\lesssim & |F^\triangle|^{1/2} k_{max} \bigg(|[K^{-1}]|\displaystyle\sum_{j=1}^3|N_{T,j}(\varphi_{0})| + \frac{1}{|\Gamma_T|^{3/2}}\nrmL[\Gamma_T]{[u_h]}\bigg)\\
 &\lesssim & h_F^{1/2} \dfrac{k_{max}}{k_{min}} \displaystyle\sum_{j=1}^3|N_{T,j}(\varphi_{0})| +k_{max} \dfrac{|F^\triangle|^{1/2}}{|\Gamma_T|^{3/2}}\nrmL[\Gamma_T]{[u_h]}.
\end{eqnarray*}  
Using the previous bounds and \eqref{eq: bound_phi_0_F} in \eqref{eq:eta_F_2}, as well as $|F^\triangle|\lesssim |\Gamma_T|$, we obtain: 
\begin{equation*}
\frac{\sqrt{h_F}}{\sqrt{k_{\Gamma}}} \nrmL[F]{(\sigIRT - K\nabla_h u_h)_{|T}\cdot n_F}\lesssim \frac{k_{max}}{k_{min}}\bigg( \frac{h_T}{\sqrt{k_{\Gamma}}} \sum_{j=1}^3|N_{T,j}(\varphi_{0})| +\frac{\sqrt{h_F}}{\sqrt{k_{\Gamma}}} \nrmL[\Gamma_T]{[K\nabla u_h\cdot n_\Gamma]} + \frac{\sqrt{h_F k_{\Gamma}}}{|\Gamma_T|}\nrmL[\Gamma_T]{[u_h]}\bigg).
\end{equation*} 
Thanks to Lemmas \ref{lem: bound of N_i,T(phi_0)},  \ref{lem: bound_2_effica} and \ref{lem: bound_3_effica}, we finally deduce that:
\begin{equation*}
\frac{\sqrt{h_F}}{\sqrt{k_{\Gamma}}} \nrmL[F]{(\sigIRT - K\nabla_h u_h)_{|T}\cdot n_F}\lesssim \max_{T'\in \Delta_T \cap \T^{\Gamma}}\frac{h_{T'}}{|\Gamma_{T'}|}\frac{k_{max}^{3/2}}{k_{min}^{3/2}}  \left(\|u-u_h\|_{h,\Delta_F}+\epsilon(\Delta_F)\right).
\end{equation*} 

We conclude by summing upon $T\in \{T_F^+,T_F^-\}$ and by using \eqref{eq:eta_F_1}.
\end{proof}

\section{Numerical simulations}\label{sec:num_sim}
To illustrate the theoretical results, we present several numerical experiments. \Blue{In \cite{Article1}, we validated both the flux implementation and the} adaptive mesh refinement (AMR) \Blue{strategy based on the error indicator $\eta_T$, which is easier to implement and computationally cheaper than $\bar \eta_T:=\eta_T+\tilde \eta_T+\sum_{F\in \F^\Gamma\cap \partial T }\eta_F$. In the following, we investigate the numerical behavior of $\eta$, $\eta_{\Gamma}$, and the error under two AMR procedures, using either $\eta_T$ or $\bar \eta_T$ as local error indicators, in order to justify the use of $\eta_T$ in \cite{Article1}. Moreover, we numerically show optimal convergence and robustness with respect to the diffusion contrast and to the mesh/interface configuration, even in situations where the theoretical constants may blow up.} We consider three test-cases.

\begin{example}[Ellipse problem]\label{example2}

Let $\Omega=[-1,1]^2$ and let $\Gamma$ be the ellipse centered at the origin, of equation $\rho=1$ where $\rho =\displaystyle{\sqrt{\frac{x^2}{a^2}+\frac{y^2}{b^2}}}$, with $2a$ the width and $2b$ the height of the ellipse. The exact solution of \eqref{eq: continuous_problem} with $g=0$ is given by
\begin{equation*}
 u(\rho)=\frac{1}{k_1}\sqrt{\rho}\quad  \text{if }\rho\leq 1,\qquad  u(\rho)=\frac{1}{k_2}\sqrt{\rho}+\frac{1}{k_1}-\frac{1}{k_2} \quad  \text{if }\rho>1.
\end{equation*}
%\begin{equation*}
%    u(x,y)=\left\{\begin{array}{ll}
%         \dfrac{1}{k_1}\rho^p&  \text{if }\rho\leq 1 \\
%        \dfrac{1}{k_2}\rho^p+\dfrac{1}{k_1}-\dfrac{1}{k_2} & \text{if }\rho>1
%    \end{array}\right..
%\end{equation*}
Here, we take $a= \displaystyle \frac{\pi}{6.18}$, $b=1.5a$, $k_1=1$ and let $k_2=\mu k_1$. For $\mu\neq 1$, the solution is merely in $H^{\frac{3}{2}-\varepsilon}(\Omega)$ for any $\varepsilon >0$ and becomes singular at the origin. 
\end{example}

We are interested in testing the influence of the additional indicator $\eta_{\Gamma}$ in the AMR procedure. We recall that $\eta_{\Gamma}^2$ is the sum of two contributions: $\sum_{F\in \F^\Gamma}\eta_F^2$, which we denote in the sequel by $\hat\eta$ and which is due to the definition of the immersed Raviart-Thomas space on the cut edges, and $\sum_{T\in\T^\Gamma}\tilde{\eta}_T^2$, which is inherent to the CutFEM formulation and which is due to the discontinuity of $u_h$ across $\Gamma$.

\begin{figure}[ht!]
    \centering  
     \begin{subfigure}[t]{0.45\textwidth}
        \centering
        \includegraphics[width=0.77\textwidth]{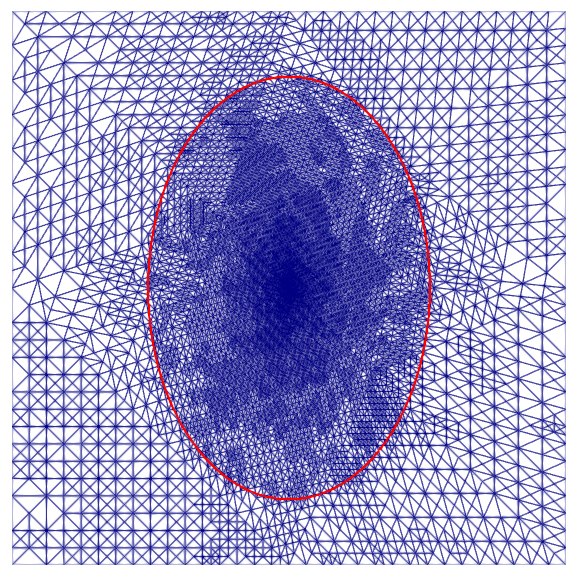}
        \caption{Final mesh (iteration 14) using $\eta_T$}
    \end{subfigure}
     \begin{subfigure}[t]{0.45\textwidth}
        \centering
        \includegraphics[width=0.77\textwidth]{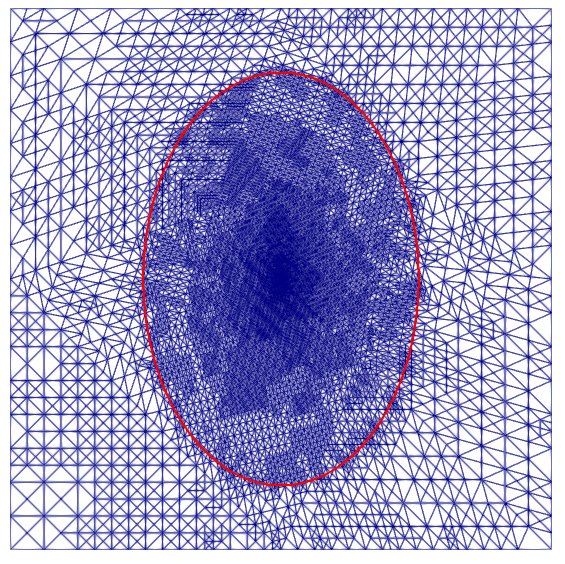}
        \caption{Final mesh (iteration 12) using $\bar \eta_T$}
    \end{subfigure}
\caption{\Cref{example2} Final adapted meshes for $\mu= 10$}
\label{fig:ex2_1}
\end{figure}

The stopping AMR criterion is that the total number of degrees of freedom $N$ is less than \Magenta{$25{,}000$}. We first take $\mu=10$ and perform two AMR procedures, one based on the main error indicator $\eta_T$ for any $T\in\T$, and the other on the indicator $\bar\eta_T=\eta_T+\tilde \eta_T+\sum_{F\in \F^\Gamma\cap \partial T }\eta_F$. Figure \ref{fig:ex2_1} shows the final meshes obtained with these refinement procedures. One can observe that they are very similar and that the refinement mainly focuses near the singularity, as expected. 
%The influence of $\eta_{\Gamma}$ seems to be negligible. 

\begin{figure}[ht!]
    \centering
 \begin{subfigure}[t]{0.49\textwidth}
        \centering
        \includegraphics[width=1\textwidth]{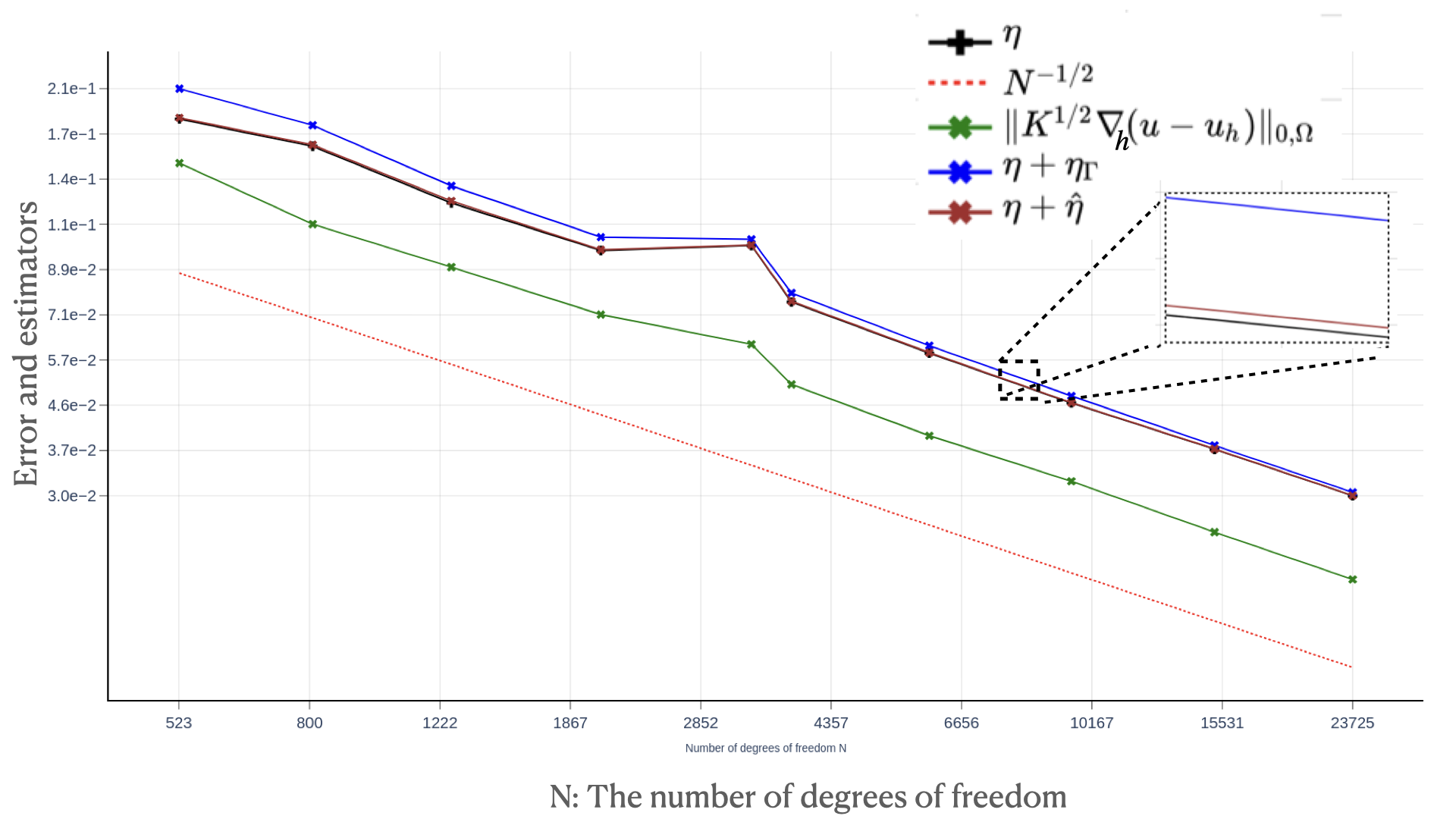}
        \caption{AMR using $\eta_T$}
    \end{subfigure}
     \begin{subfigure}[t]{0.49\textwidth}
        \centering
        \includegraphics[width=1\textwidth]{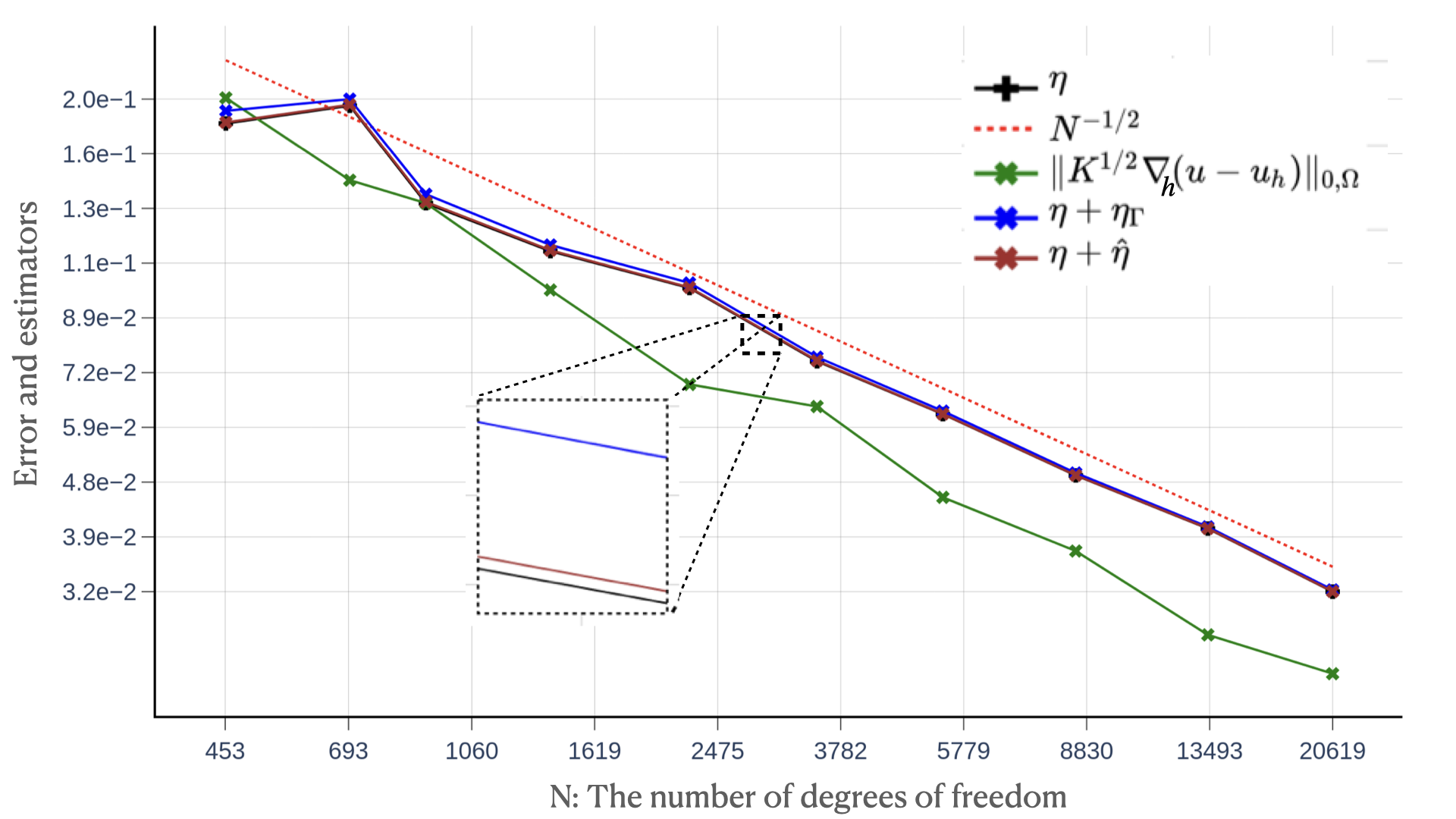}
        \caption{AMR using $\bar \eta_T$}
    \end{subfigure}
\caption{\Cref{example2} Convergence rate of \Magenta{the error}, $\eta$, $\eta+\hat \eta$ and $\eta+\eta_{\Gamma}$ for $\mu=10$}
\label{fig:ex2_2}
\end{figure}

We also compare the estimators $\eta$, $\eta+\hat \eta$ and $\eta+\eta_\Gamma$. Their decay rate for the last ten refinement iterations is shown in Figure \ref{fig:ex2_2}, while in Table \ref{tab:ex2} we give the values of $\eta$, $\eta+\eta_\Gamma$ \Magenta{and the effectivity index} for the last iterations. One can see that the two AMR procedures yield similar results and moreover, the additional estimator $\eta_{\Gamma }$ is quite small. We conclude that $\eta_{\Gamma}$ can be neglected in the AMR, which substantially simplifies the latter. \Magenta{The effectivity index, i.e. the ratio between the employed global indicator and the $H^1$ semi-norm of the error  $err=\|K^{1/2}\nabla_h (u-u_h)\|_\Omega$, is slightly smaller when using the indicator $\tilde \eta_T$ in the AMR. This is in agreement with the theoretical results, since the reliability was established for $\eta+\eta_{\Gamma}$.}

\begin{figure}[ht!]
    \centering  
    \begin{subfigure}[t]{0.65\textwidth}
        \centering
        \includegraphics[width=0.85\textwidth]{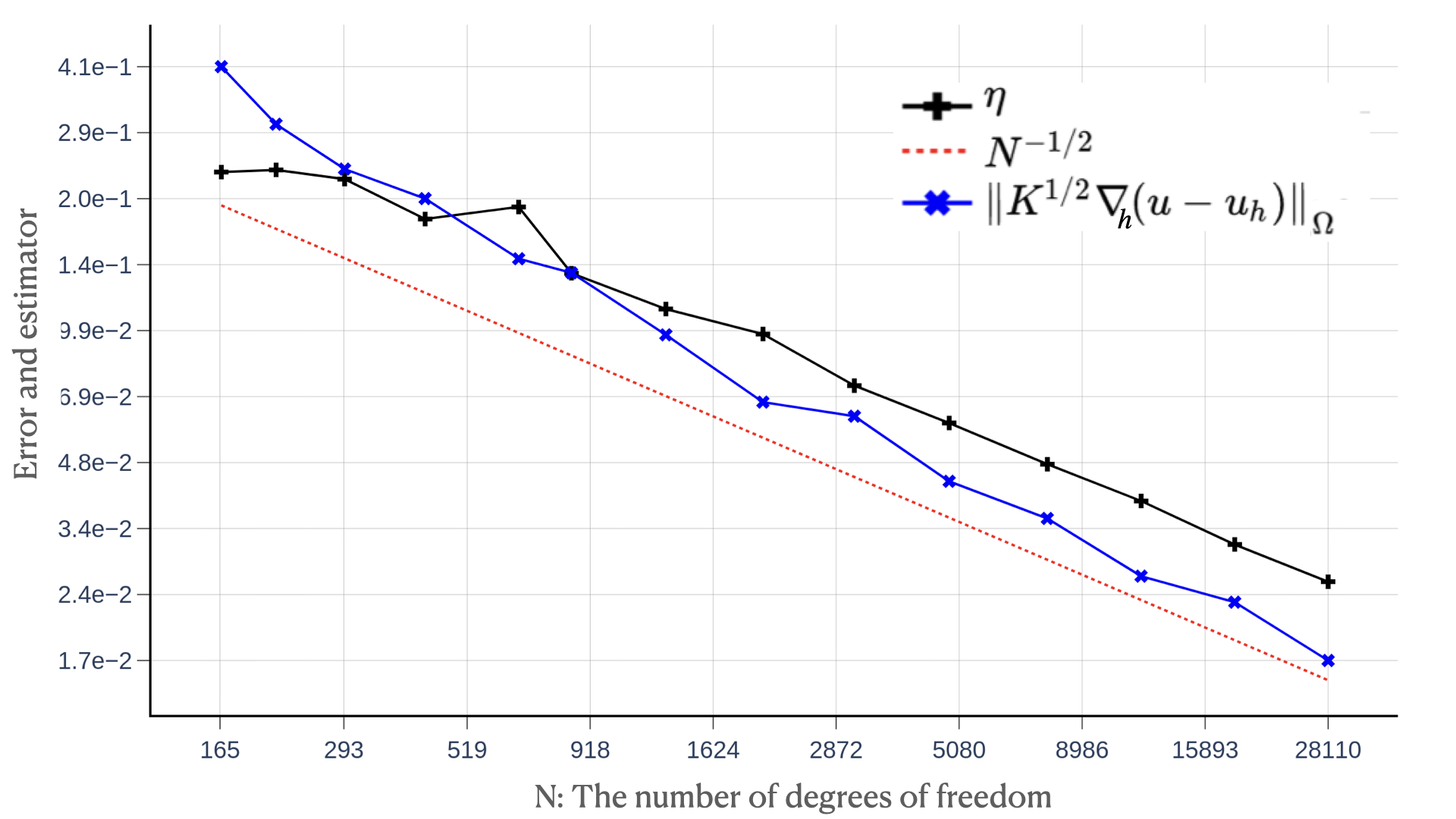}
        \caption{Convergence rate for the exact error and $\eta$}
    \end{subfigure}%
    \begin{subfigure}[t]{0.35\textwidth}
        \includegraphics[width=0.95\textwidth]{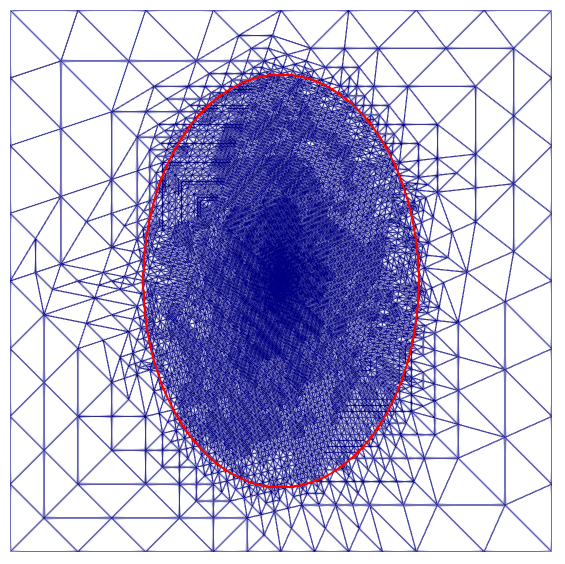}
        \caption{Final mesh (iteration 13)}
    \end{subfigure}
\caption{\Cref{example2} AMR using $\eta_T$: \Red{estimator} and error convergence \Red{(left)} and final mesh \Red{(right)} for $\mu=10^6$ }
\label{fig:ex2_3}
\end{figure}

Finally, we consider a large jump in the diffusion coefficients, that is $\mu=10^6$. Figure \ref{fig:ex2_3} shows the convergence rate for both the error and the estimator $\eta$, as well as the final mesh obtained using only  the indicator $\eta_T$ for the adaptive refinement. Although the ratio $k_{max}/k_{min}$ is large, and hence the (theoretical) efficiency constant $C_T$ is large (independently of the mesh/interface cut), this does not seem to influence the estimator behavior since we retrieve the optimal convergence rate $O(N^{-1/2})$.

\begin{table}[ht!]
\centering
\begin{minipage}{0.48\textwidth}
\centering
\begin{tabular}{|c|c|c|c|}
\hline
{$N$} & {$\eta+\eta_\Gamma$} &  {$\eta$} & \Magenta{$\eta/err$}\\
\hline
3830  &$ 7.9\times 10^{-2}$ &$ 7.5\times 10^{-2}$&$\Magenta{ 1.48}$ \\
5995  &$ 6.1\times 10^{-2}$ &$ 5.9\times 10^{-2}$&$ \Magenta{1.4}$ \\
9511  &$ 4.8\times 10^{-2}$ &$ 4.6\times 10^{-2}$&$ \Magenta{1.45}$ \\
15142 &$ 3.8\times 10^{-2}$ &$ 3.7\times 10^{-2}$&$ \Magenta{1.49}$ \\
23725 & $3\times10^{-2}$ & $2.9\times10^{-2}$ &$ \Magenta{1.49}$\\
\hline
\end{tabular}
\caption*{(a) AMR using only $\eta_T$}
\end{minipage}
\hfill 
\begin{minipage}{0.48\textwidth}
\centering
\begin{tabular}{|c|c|c|c|}
\hline
{$N$} & {$\eta+\eta_\Gamma$} & {$\eta$} & \Magenta{$(\eta+\eta_\Gamma)/ err$}\\
\hline
3486  &$7.6\times 10^{-2}$ &  $7.5\times 10^{-2}$&$ \Magenta{1.20}$ \\
5379  &$6.2\times 10^{-2}$ &  $6.1\times 10^{-2}$&$ \Magenta{1.38}$ \\
8506  &$4.9\times 10^{-2}$ &  $4.9\times 10^{-2}$&$ \Magenta{1.34}$\\
13410 &$4.0\times 10^{-2}$ &  $4.0\times 10^{-2}$&$ \Magenta{1.50}$ \\
20619 &$3.1\times 10^{-2}$ &  $3.1\times 10^{-2}$ &$ \Magenta{1.37}$ \\
\hline
\end{tabular}
\caption*{(b) AMR using $\bar \eta_T$}
\end{minipage}
\caption{\Cref{example2} Estimators \Magenta{and effectivity index}: last iterations for $\mu=10$}
\label{tab:ex2}
\end{table}

\begin{example} [Sinusoidal interface problem]\label{example3}
We next consider an interface problem in the square domain \( \Omega = [-1, 1]^2 \), with several interfaces  described by the periodic level set function $\phi(x, y) = \sin(2\pi x) \cos(2\pi y) - 0.2$. Thus, \( \Gamma = \{ (x, y) \in \Omega;\, \phi(x, y) = 0 \} \). The diffusion coefficient is:
\[
k(x, y) =1\quad \text{if }\phi(x, y) < 0,\qquad k(x, y) =\Blue{100}\quad \text{if }\phi(x, y) \geq 0.
\]
The exact solution is given by $u(x, y) =\displaystyle\frac{1}{k(x,y)} \phi(x, y)$ for $(x, y) \in \Omega$. %The solution $u$ and the source term $f=-\Delta\phi $ are shown in Figure \ref{fig:ex3_1}.
\end{example}

\begin{figure}[ht!]
    \centering
    \begin{subfigure}[t]{0.6\textwidth}
        \includegraphics[width=0.9\textwidth]{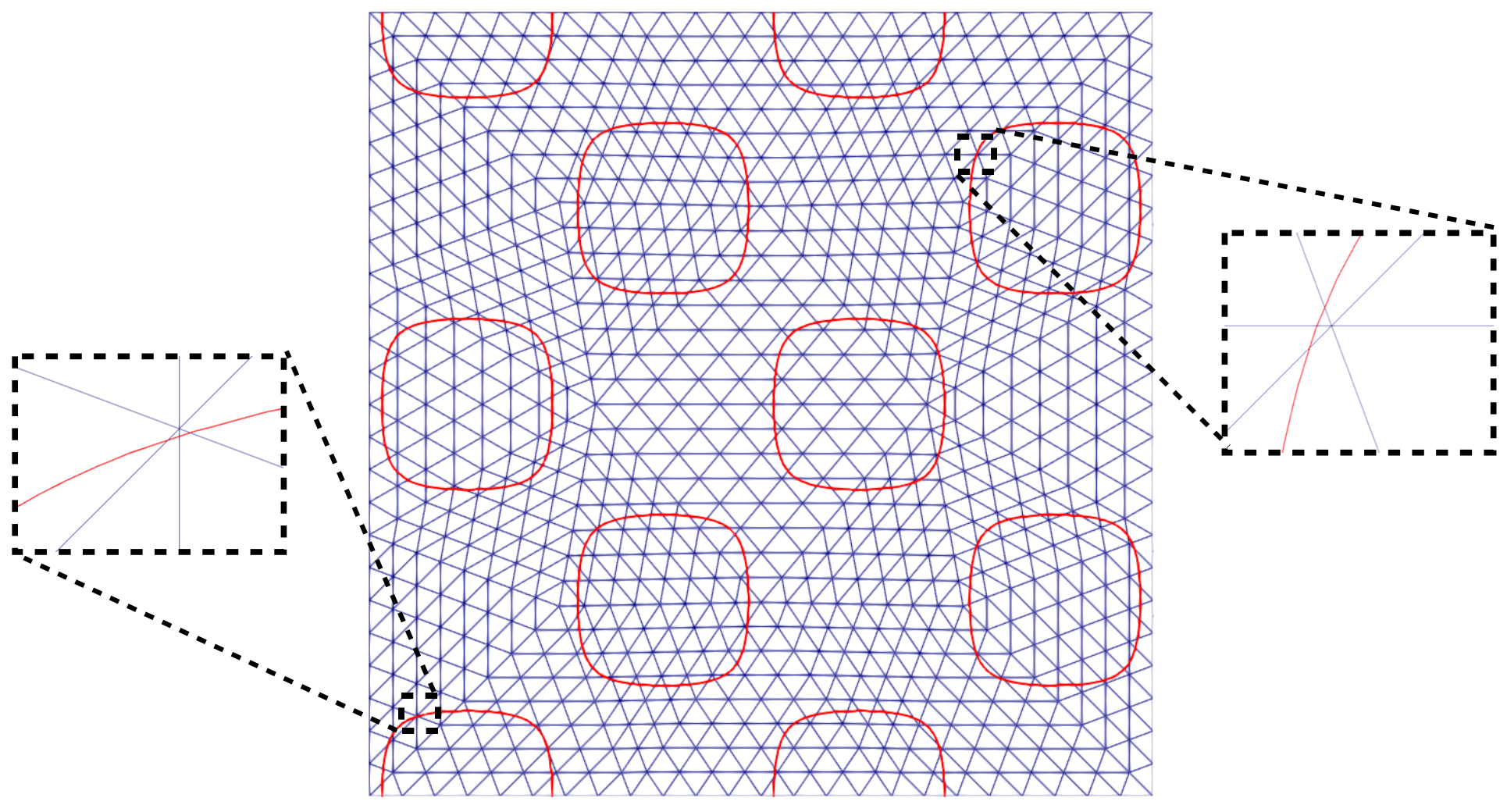}
        \caption{Initial mesh (iteration 0)}
    \end{subfigure}%
     \begin{subfigure}[t]{0.33\textwidth}
        \centering
        \includegraphics[width=0.9\textwidth]{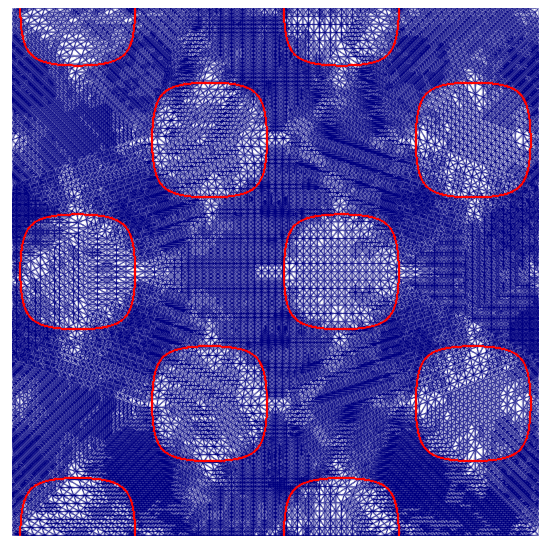}
        \caption{ Final mesh (iteration \Blue{13})}
    \end{subfigure}
\caption{\Cref{example3} Initial and final meshes}    
\label{fig:ex3_2}     
\end{figure}

\begin{figure}[ht!]
\centering
\includegraphics[width=.6\textwidth]{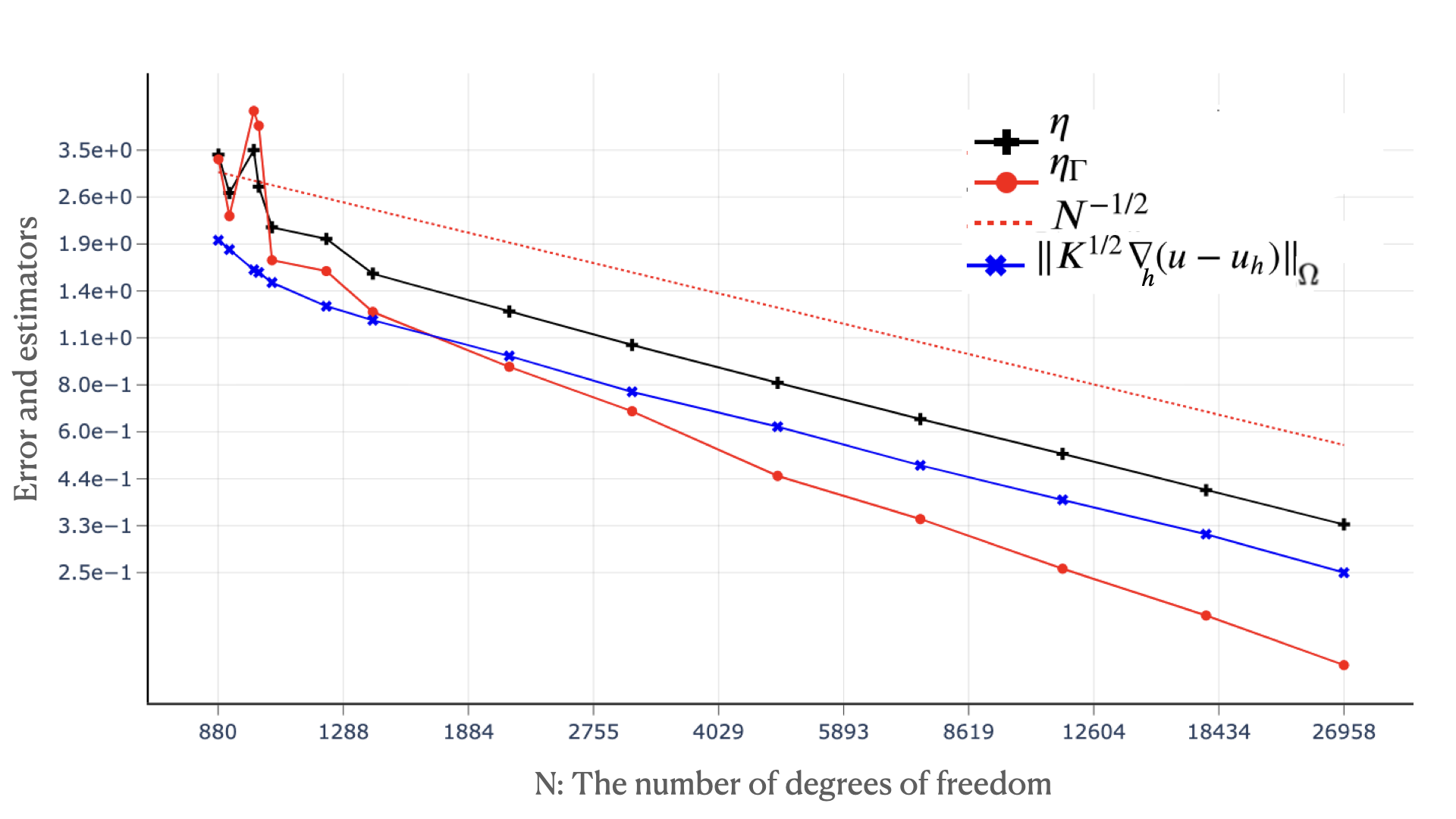}
²\caption{\Cref{example3} Convergence rate for the error and the estimators $\eta$ and $\eta_{\Gamma}$}
\label{fig:ex3_3} 
\end{figure}

\begin{figure}[ht!]
\centering
\includegraphics[width=0.86\textwidth]{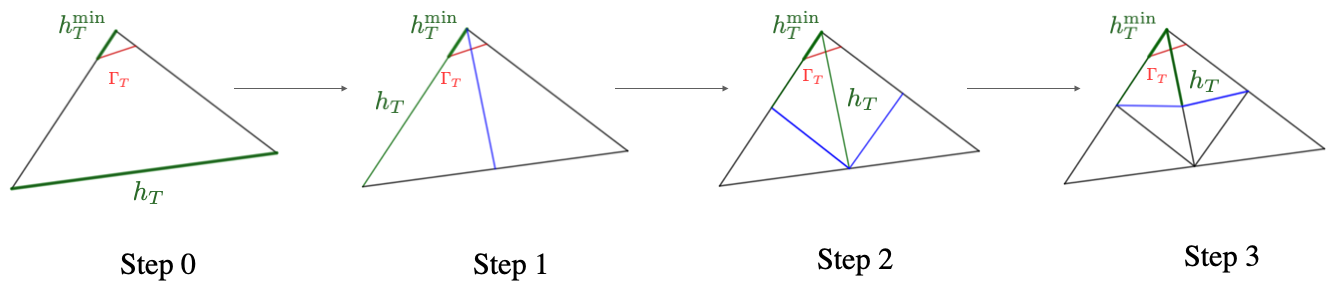}
\caption{Illustration of the evolution of $\dfrac{h_T}{\sqrt{h_T^{min}|\Gamma_T|}}$ during the refinement of a cut cell $T$}
\label{fig:cut_setps1} 
\end{figure}

%This example is particularly suitable for testing smooth but oscillatory interface geometries. 
In this example, we use only $\eta_T$ as error indicator in AMR. \Magenta{The stopping criterion is that $N\le 30{,}000$.} In Figure \ref{fig:ex3_2} (a) we show the initial mesh with the periodic interface, as well as two zooms near the interface, while in Figure \ref{fig:ex3_2} (b) the final mesh is given. As can be seen in the zooms, there are multiple bad cuts which yield large ratios $h_T/\sqrt{h_T^{min}|\Gamma_T|}$, and consequently large theoretical efficiency constants $C_T$ on the affected triangles. \Blue{In addition, the multiplicative factor $(k_{max}/k_{min})^{3/2}$ in the efficiency constants $C_T$ and $C_F$ is equal to $10^3$. This explains the peak observed at iterations 2 and 3 in Figure \ref{fig:ex3_3} for the estimators $\eta$, and particularly $\eta_{\Gamma}$, which also contains the multiplicative factor $\sqrt{h_T}/\sqrt{h_T^{min}}$. It is important to note that these ratios decrease during the AMR procedure and eventually become bounded. Figure \ref{fig:cut_setps1} illustrates several refinement steps of a cut triangle (in the worst case configuration) using the longest-edge refinement method, showing how $h_T/\sqrt{h_T^{min}|\Gamma_T|}$ diminishes after a few steps. Although such situations may arise in the early iterations, the estimator is able to cope with them and its overall performance is not affected.} As seen in Figure \ref{fig:ex3_3}, we recover the optimal convergence rate $O(N^{-1/2})$ for both the error and the estimator $\eta$.

%\begin{figure}[ht!]
%    \centering
%    \begin{subfigure}[t]{0.55\textwidth}
%        \includegraphics[width=1\textwidth]{New-Example/U_exact.png}
%        \caption{Analytical solution $u$}
%    \end{subfigure}%
%    \begin{subfigure}[t]{0.45\textwidth}
%        \centering
%        \includegraphics[width=1\textwidth]{New-Example/f.png}
%        \caption{Source term $f$}
%    \end{subfigure}
%\caption{\Cref{example3} Exact solution and data}    
%\label{fig:ex3_1}   
%\end{figure}  

\begin{figure}[ht!]
    \centering
    \begin{subfigure}[t]{0.35\textwidth}
        \includegraphics[width=0.88\textwidth]{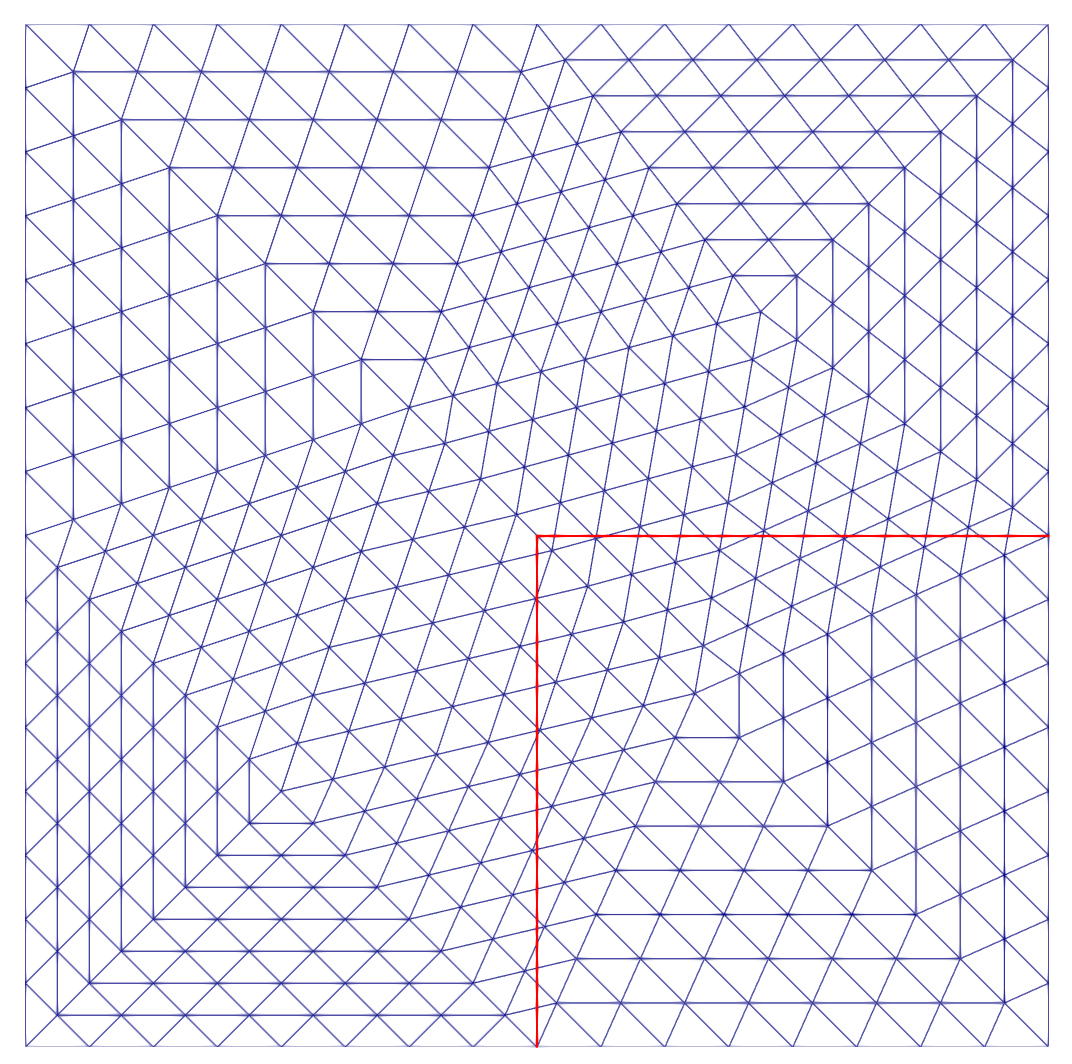}
        \caption{Initial mesh (iteration 0)}
    \end{subfigure}\quad 
     \begin{subfigure}[t]{0.35\textwidth}
        \centering
        \includegraphics[width=0.88\textwidth]{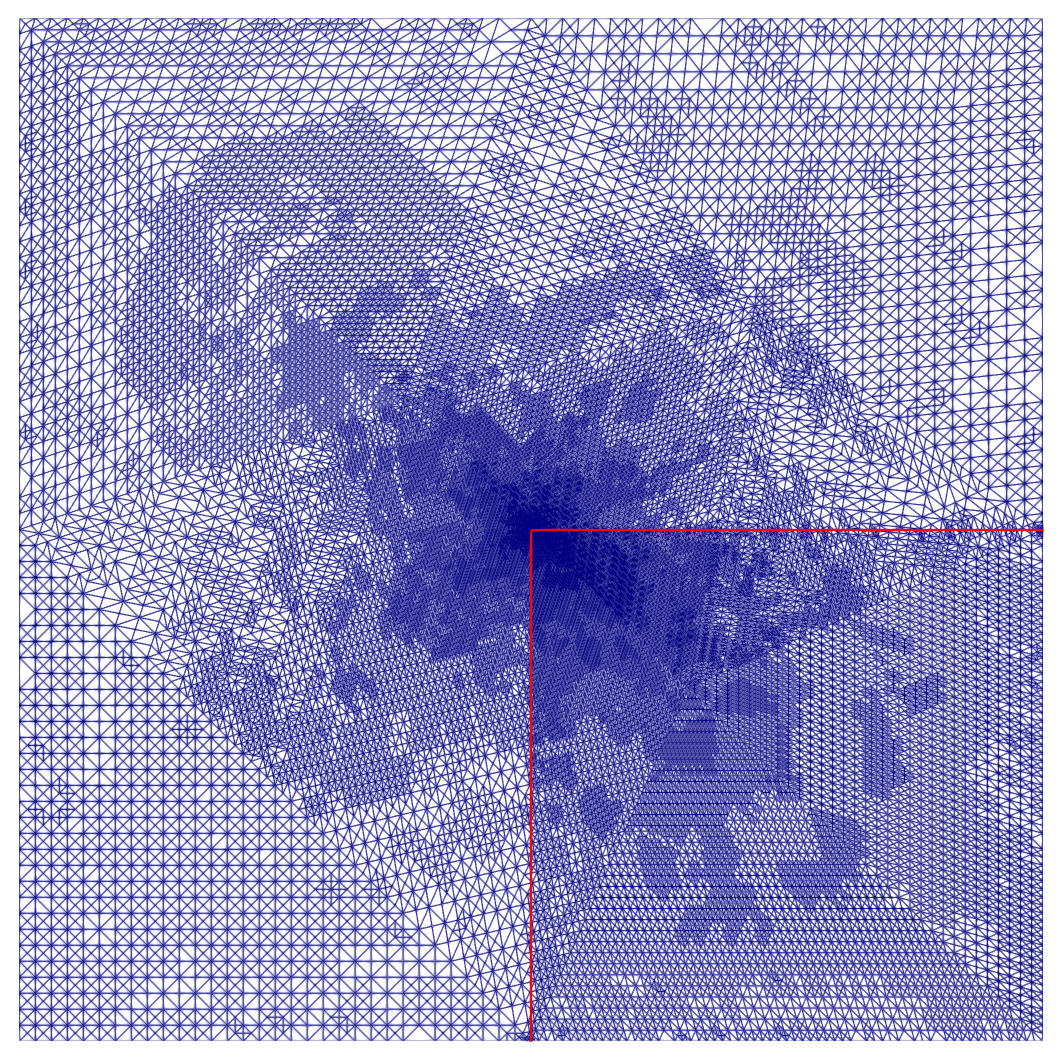}
        \caption{Final mesh (iteration 10)}
    \end{subfigure}
\caption{\Cref{example4} Initial and final meshes}    
\label{fig:ex4_2}     
\end{figure}

\begin{figure}[ht!]
\centering
\includegraphics[width=0.65\textwidth]{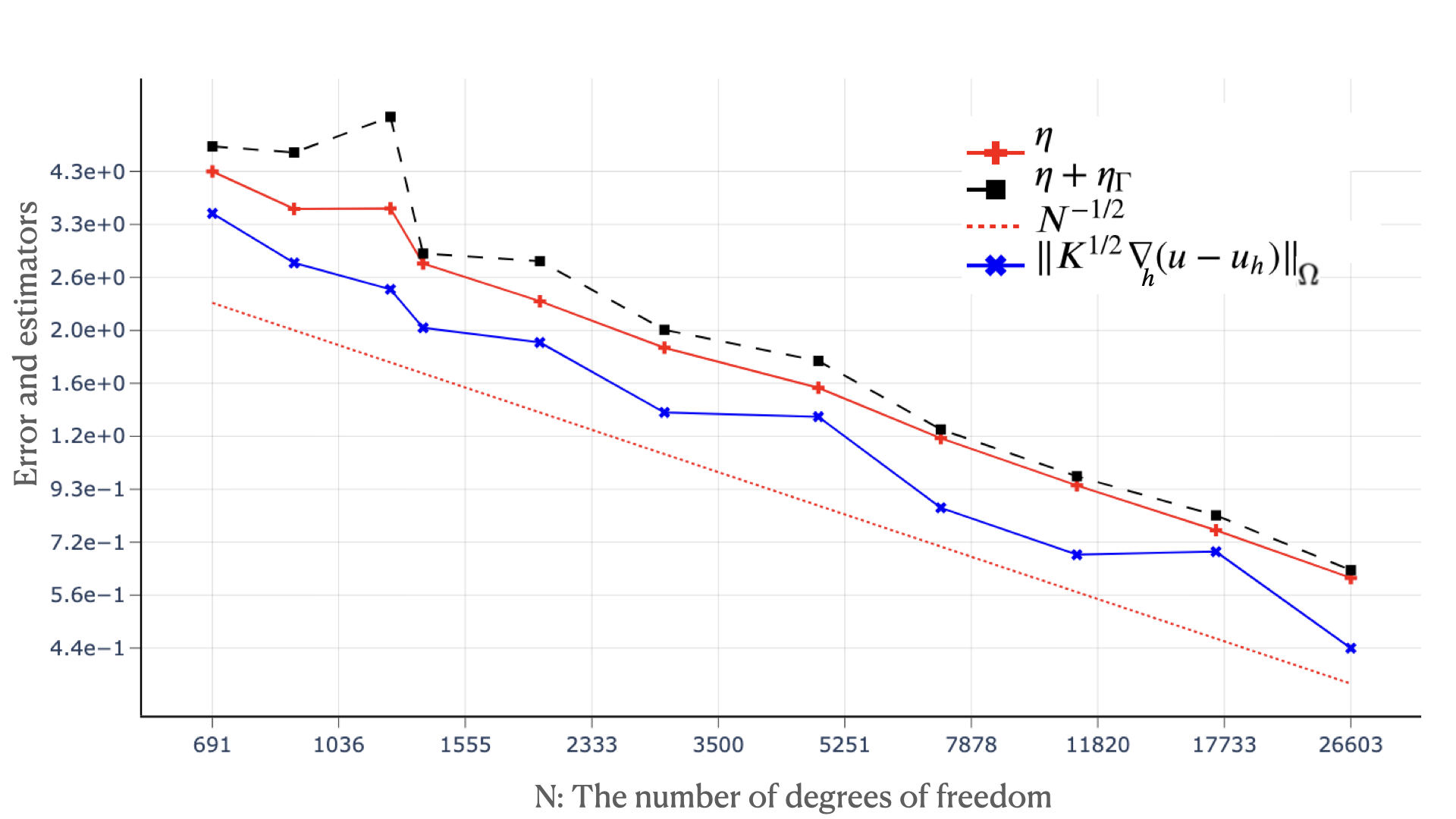}
\caption{\Cref{example4} Convergence rate for the error and the estimators $\eta$ and $\eta+\eta_{\Gamma}$}
\label{fig:ex4_3} 
\end{figure}

\Blue{\begin{example}[Non-smooth interface problem]\label{example4}
Finally, we consider a non-smooth interface $\Gamma$ which presents a corner. Let $\Omega = [-1,1]^2$ and $\Gamma = \{(x,0);\, x \in [0,1]\}\cup\{(0,y);\, y \in [-1,0]\}$. The interface can be described by the following level-set function:
$\phi(x,y)=y$ if $y > -x$, $\phi(x,y)=-x$ otherwise. The exact solution is given by
\[
u(x,y)=
\begin{cases}
\dfrac{25}{k_1}\, x y (x^2+y^2)^{3/4} & \text{if } \phi(x,y) < 0,\\[6pt]
\dfrac{25}{k_2}\, x y (x^2+y^2)^{3/4} & \text{if } \phi(x,y) \ge 0,
\end{cases}
\]
and satisfies homogeneous transmission conditions and non-homogeneous Dirichlet boundary condition.
\end{example}
We set $k_1 = 1$ and $k_2 = 10$. Note that the solution has a singularity at the origin, which coincides with the corner of the interface. 
The adaptive refinement is performed using the indicator $\eta_T$ and the stopping criterion is $N\le 30{,}000$. Figures \ref{fig:ex4_2}(a) and \ref{fig:ex4_2}(b) show the initial and final meshes. One can see that the refinement occurs mainly near the origin, as expected. In Figure \ref{fig:ex4_3}, we observe a behavior similar to that in Figure \ref{fig:ex3_3}: the estimators $\eta$, and particularly $\eta_{\Gamma}$, exhibit a peak at refinement iteration 2, due to some bad cuts by the interface. Nevertheless, the optimal convergence rate $O(N^{-1/2})$ is recovered.
}


\textbf{Acknowledgement.} This project has received funding from the European Union’s Horizon H2020 Research and Innovation under Marie Curie Grant Agreement N° 945416.

%%-----------------------------

%%-----------------------------
\end{document}